\documentclass[12pt,a4paper]{article}
 \usepackage{emlines2,amstext,amssymb,amscd}

\begin{document}

\title{Arrays and the octahedron recurrence}
\author{V. I. Danilov and G.A. Koshevoy \thanks{We thank A.Vershik and
I.Pak for useful discussions
and comments. A partial support from the grant NSch-1939.2003.6 is
acknowledged. G.A. Koshevoy also thanks for support the Foundation
of Support of Russian Science. }}

\date{}

\maketitle

\section{Introduction}

Recently, in  \cite{oct,h-k,PV} several interesting bijections
have been constructed. In \cite{oct,h-k} bijections relate special
sets of discretely concave functions (hives) on triangular grids
and the octahedron recurrence (OR) plays the main role for these
bijections. Bijections in \cite{PV} relate special sets of Young
tableaux and constructions of these bijections based on standard
algorithms in this theory, {\em jeu de taquen}, Schutzenberger
involution, {\em tableaux switching} etc.

The transparentness of the octahedron recurrence is a undoubtable
advantage, but a rationale why the OR provides natural bijections
(and even bijections) was rather obscured.

In this paper we investigate these constructions from the third
point of view, combinatorics of arrays, theory worked out by the
authors in \cite{umn}. Arrays naturally related as well to
functions on the lattice of integers as to Young tableaux, and
have some advantages comparing to functions and tableaux. For
example, Young tableaux are nothing but integer-valued  {\bf
D}-tight arrays. In the tensor category of arrays, the bijections
of associativity and commutativity arise naturally. We establish
coincidence of these bijections with that defined in
\cite{oct,h-k,PV}.

In order to relate different approaches and to reveal
combinatorics of the octahedron recurrence, we, first, show that
the octahedron recurrence agrees with discrete convexity and,
second, we construct another bijection using the OR, the
functional form of the RSK correspondence.

The paper is organized as follows: after a brief introduction to
the octahedron recurrence and discrete concave functions, in
Section 4 we state Theorem 1 on heredity of discrete concavity
under propagation due to octahedron recurrence. In Section 5 we
recall definitions and facts from theory of arrays, which are of
use in this paper. In Theorem 2 we establish a relation between
the OR and the operation of condensation of arrays. This might be
seen as a functional form of the RSK correspondence (more
precisely, modified RSK, see \cite{umn}). In Section 7 we recall
the definitions of the associativity bijection for arrays and the
associativity bijection from \cite{oct}. In Theorem 3 (Section 8)
we state that these bijections coincide. More subtle constructions
are used for proving Theorem 4, in which we establish coincidence
of our bijection of commutativity with the functional
commutativity bijection in \cite{h-k} and with two fundamental
symmetries of Pak and Vallexo (\cite{PV}, Conjecture 1).

\section{Octahedron recurrence}

The main idea of the octahedron recurrence is rather transparent.
Specifically, consider the octahedron

\unitlength=.7mm \special{em:linewidth 0.4pt} \linethickness{0.4pt}
\begin{picture}(83.00,64.00)(-30,-5)
\emline{65.00}{20.00}{1}{40.00}{20.00}{2}
\emline{40.00}{20.00}{3}{55.00}{35.00}{4}
\emline{55.00}{35.00}{5}{80.00}{35.00}{6}
\emline{80.00}{35.00}{7}{65.00}{20.00}{8}
\emline{65.00}{20.00}{9}{60.00}{50.00}{10}
\emline{60.00}{50.00}{11}{60.00}{50.00}{12}
\emline{60.00}{50.00}{13}{55.00}{35.00}{14}
\emline{55.00}{35.00}{15}{60.00}{5.00}{16}
\emline{60.00}{5.00}{17}{65.00}{20.00}{18}
\emline{80.00}{35.00}{19}{60.00}{49.00}{20}
\emline{60.00}{49.00}{21}{40.00}{20.00}{22}
\emline{40.00}{20.00}{23}{60.00}{5.00}{24}
\emline{60.00}{5.00}{25}{80.00}{35.00}{26}
\put(83.00,35.00){\makebox(0,0)[cc]{$a'$}}
\put(60.00,52.00){\makebox(0,0)[cc]{{\bf 1}}}
\put(49.00,38.00){\makebox(0,0)[cc]{$b'$}}
\put(37.00,20.00){\makebox(0,0)[cc]{$a$}}
\put(72.00,17.00){\makebox(0,0)[cc]{$b$}}
\put(60.00,2.00){\makebox(0,0)[cc]{{\bf 0}}}
\end{picture}

\hfill Picture 1. \hfill \medskip

\noindent with the vertexes ${\bf 0}, a,a',b,b'$ and  {\bf 1}. Let
$f$ be a real-valued function given at the points ${\bf 0},
a,a',b,b'$. Then we can {\em propagate}  $f$ to the point {\bf 1}
by the following rule
$$
f({\bf 1})=\max(f(a)+f(a'), f(b)+f(b'))-f({\bf 0}).
$$
We refer to  \cite{Sp} for justification of this rule and its
interesting appearances in combinatorics. Rather unexpectedly this
related to flips in \cite{fock}. We want to point out a relation
of this rule to concavity. Specifically, suppose
$f(b)+f(b')=\max(f(a)+f(a'), f(b)+f(b'))$. Then, we have $f({\bf
0})+f({\bf 1})=f(b)+f(b')$. This means that the restriction of the
function to the rhombus ${\bf 0},b,{\bf 1},b'$ coincides with the
restriction of an affine function $h$. Moreover,
$$
h(a)+h(a')=2h((a+a')/2)=2h((b+b')/2)= $$
$$
2f((b+b')/2)=f(b)+f(b')\ge f(a)+ f(a').
$$
We can choose  $h$ such that  $f(a)\le h(a)$ and $f(a')\le h(a')$
hold true. That means that the function  $f$ is sub-affine on the
octahedron, i.e. $f$ looks alike a {\em concave} function.
Moreover, the rhombus  ${\bf 0},b,{\bf 1},b'$ is an affinity set
of $f$.

In other words, we propagate the function  $f$ to the point  {\bf
1} in order to get a concave (discretely) function on the
octahedron, such that an affinity area (the convex hull of the
affinity set) has to contain the vector ${\bf 0}{\bf 1}$, a {\em
propagation vector}.

Now, using this rule, which is called the {\em octahedron
recurrence (OR)}, we can propagate a function given at some domain
to a large domain. Here is one of possible initial domains (see
\cite{Sp}). Let us consider the set $L$ of points $(n,i,j)$ with
integers $i,j,n$, $n\ge 0$ and $n=i+j$ ($\text{mod} \ 2$). Suppose
a function  $f$ is given at a subset of  $L$ constituted from
points of  the form  $(\cdot ,\cdot ,0 )$ and $(\cdot ,\cdot ,1)$.
Then using the octahedron recurrence with the propagation vector
$(0,0,2)$ we can propagate the function to points of  $L$ of the
form  $(\cdot ,\cdot , 2 )$, on the next step to points of the
form  $(\cdot ,\cdot , 3 )$ and so on to the whole  $L$.

Of course, the initial data can be given at more sophisticated
subsets, see  \cite{Sp} and \cite{h-k}.\medskip

We will display the stuff in a slightly different manner
\footnote{Sometimes it is convenient to draw pictures for the OR
with different propagation vectors. In order to set the octahedron
recurrence, we have to choose a unimodular set in the lattice
$\mathbb Z^3$, say, $\{e_1$, $e_2$, $e_3$, $e_1-e_2$, $e_1-e_3$,
$e_2-e_3\}$ and the propagation vector  $e_3-e_1+e_2$, where
$e_1$, $e_2$, $e_3$ is a basis in the lattice $\mathbb Z^3$. Then
the primitive octahedron becomes the convex hull of the points
$0$, $e_3-e_1$, $e_2$, $e_3$, $e_2-e_1$, $e_3-e_1+e_2$. An integer
translation of a plane, spanned by a triple of vectors in the
unimodular set, is a {\em modular flat}. A {\em non-modular flats}
are parallel to planes spanned either by the pair ($e_3-e_1$,
$e_2$) or ($e_2-e_1$, $e_3$).}. Specifically, we consider the
integer orhtant $\mathbb Z^3_+$ (with coordinates  $x,y,z$). The
propagation vector is  $(1,0,1)$, that is proportional to the
vector  $OD$ see Picture 2. On Picture 2 with $n=m$, we can see
modular and non-modular flats: the modular flats are parallel to
the faces of the tetrahedron $OEAB$, and the non-modular flats
parallel to the face $OEDC$ and the plane passing through $ODB$.

We locate the initial data of functions on  $Oxz$, $Oyz$
(typically equal $0$ on $OEA$ and $OEDC$) and on the plane $z=y$
(more precisely at integer points of the rectangle  $OABC$, and
here are the main data).

\unitlength=.700mm \special{em:linewidth 0.4pt}
\linethickness{0.4pt}
\begin{picture}(93.00,70)(-30,0)
\put(30.00,20.00){\vector(3,1){60.00}}
\put(30.00,20.00){\vector(0,1){43.00}}
\put(30.00,20.00){\vector(3,-1){50.00}}
\put(33.00,64.00){\makebox(0,0)[cc]{$z$}}
\put(85.00,3.00){\makebox(0,0)[cc]{$x$}}
\put(90.00,36.00){\makebox(0,0)[cc]{$y$}}
\emline{30.00}{50.00}{1}{60.00}{60.00}{2}
\emline{60.00}{60.00}{3}{90.00}{50.00}{4}
\emline{60.00}{40.00}{5}{90.00}{50.00}{6}
\emline{60.00}{40.00}{7}{30.00}{50.00}{8}
\put(26.00,17.00){\makebox(0,0)[cc]{$O$}}
\put(26.00,50.00){\makebox(0,0)[cc]{$E$}}
\put(62.00,63.00){\makebox(0,0)[cc]{$A$}}
\put(106.00,52.00){\makebox(0,0)[cc]{$B=(n,m,m)$}}
\put(60.00,43.00){\makebox(0,0)[cc]{$D$}}
\put(58.00,7.00){\makebox(0,0)[cc]{$C$}}
\emline{31.00}{50.00}{9}{88.00}{50.00}{10}
\emline{60.00}{40.00}{11}{60.00}{10.00}{12}
\emline{60.00}{10.00}{13}{90.00}{50.00}{14}
\put(30.00,20.00){\vector(3,2){29.00}}
\bezier{20}(60.00,60.00)(46.00,42.00)(30.00,20.00)
\bezier{20}(30.00,20.00)(60.00,35.00)(90.00,50.00)
\end{picture}

\hfill Picture 2. \hfill \medskip

Due to the octahedron recurrence propagation we get a function on
the prism $OEABDC$, and of our particular interest will be the
resulting functions at the rectangle $EABD$ and at triangle $BCD$.
For $n=m$, we will be also interested for functions at the
tetrahedron $OBAE$ and the half-octahedron  $OEDCB$.

In this set-up, the unit octahedron if of the form depicted at
Picture 3.

\unitlength=.7mm \special{em:linewidth 0.4pt} \linethickness{0.4pt}
\begin{picture}(83.00,45.00)(-40,0)
\emline{70.00}{10.00}{1}{50.00}{15.00}{2}
\emline{50.00}{15.00}{3}{50.00}{30.00}{4}
\emline{50.00}{30.00}{5}{70.00}{25.00}{6}
\emline{70.00}{25.00}{7}{70.00}{10.00}{8}
\emline{50.00}{15.00}{9}{40.00}{5.00}{10}
\emline{40.00}{5.00}{11}{50.00}{30.00}{12}
\emline{50.00}{30.00}{13}{80.00}{35.00}{14}
\emline{80.00}{35.00}{15}{70.00}{25.00}{16}
\emline{70.00}{25.00}{17}{40.00}{5.00}{18}
\emline{40.00}{5.00}{19}{70.00}{10.00}{20}
\emline{70.00}{10.00}{21}{80.00}{35.00}{22}
\emline{80.00}{35.00}{23}{50.00}{15.00}{24}
\put(37.00,3.00){\makebox(0,0)[cc]{$a$}}
\put(83.00,37.00){\makebox(0,0)[cc]{$a'$}}
\put(48.00,15.00){\makebox(0,0)[cc]{{\bf 0}}}
\put(73.00,24.00){\makebox(0,0)[cc]{{\bf 1}}}
\put(48.00,32.00){\makebox(0,0)[cc]{$b'$}}
\put(73.00,8.00){\makebox(0,0)[cc]{$b$}}
\put(50.00,15.00){\vector(2,1){18.00}}
\end{picture}

\hfill Picture 3. \hfill \medskip

Thus, a primitive propagation takes the following form: given
values at the points  ${\bf 0}$, $a$ and $b$ at the ground flour
and two values at the points $a'$ and $b'$ at the first flour, due
to the OR we get a value at the third point {\bf 1} at the first
flour. If points  ${\bf 0}$, $b$, $b'$ and {\bf 1} are located at
the qudrant $Oxz$ (that is they have the $y$-th coordinate equals
$0$), we set the value in  {\bf 1} by the rule: $f({\bf
1})=f(b)+f(b')-f({\bf 0})$ (an instance of the octahedron
recurrence for the case $f=-\infty$ for points outside the orhtant
$\mathbb Z^3_+$). \medskip

We claim that functions, which we get as an output of the
octahedron recurrence, inherit  some concavity properties of input
functions. The next two sections are devoted to this issue.

\section{Discrete concave functions on  {\bf 2D}-grids}

We consider functions on $\mathbb Z^2$ defined on finite sets of
special form. We call such sets { grids} and they are specified as
follows. A finite subset $T\subset \mathbb Z^2$ is a {\em grid} if
i) $T$ has no holes, i.e. $T=co(T)\cap \mathbb Z^2$, and ii) any
edge of the convex hull $co(T)$ is parallel to one of the vectors
$(1,0)$, $(0,1)$ ш $(1,1)$. (Obviously, a grid has a hexagonal
shape, which might degenerated to a pentagon, a trapezoid, a
parallelogram or a triangle.)

Let  $f:T \to \mathbb R $ be a function on a grid $T$. A primitive
triangle in $T$ is either a triple  $x$, $x+(0,1)$ and $x+(1,1)$
of points of $T$, or a triple  $x$, $x+(1,0)$, $x+(1,1)$. Convex
hulls of these primitive triangles constitute a simplicial
decomposition of  co$(T)$ (if $T$ is not one-dimensional). We
uniquely interpolate the function  $f$ by affinity to the
triangles on this decomposition of co$T$, and get a function
$\tilde f: co(X)\to\mathbb R$.\medskip

{\bf  Definition}. A function  $f$ on a grid  $T$ is said to be
{\em discrete concave}, if the interpolation $\tilde f$ is a
concave function on co$(T)$.\medskip

We can reformulate discrete concavity of a function $f$ without
using the interpolation $\tilde f$. Namely we have to require
validity of three types of ``rhombus'' inequalities. Consider
``primitive'' rhombus in $T$ of the form

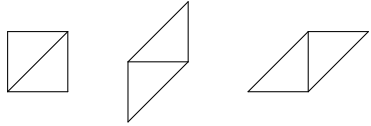
\unitlength=.8mm \special{em:linewidth 0.4pt} \linethickness{0.4pt}
\begin{picture}(100.00,30.00)(-10,0)
\emline{60.00}{5.00}{1}{60.00}{15.00}{2}
\emline{60.00}{15.00}{3}{70.00}{25.00}{4}
\emline{70.00}{25.00}{5}{70.00}{15.00}{6}
\emline{70.00}{15.00}{7}{60.00}{5.00}{8}
\emline{60.00}{15.00}{9}{70.00}{15.00}{10}
\emline{50.00}{10.00}{11}{40.00}{10.00}{12}
\emline{40.00}{10.00}{13}{40.00}{20.00}{14}
\emline{40.00}{20.00}{15}{50.00}{20.00}{16}
\emline{50.00}{20.00}{17}{50.00}{10.00}{18}
\emline{50.00}{20.00}{19}{40.00}{10.00}{20}
\emline{90.00}{10.00}{21}{100.00}{20.00}{22}
\emline{100.00}{20.00}{23}{90.00}{20.00}{24}
\emline{90.00}{20.00}{25}{80.00}{10.00}{26}
\emline{80.00}{10.00}{27}{90.00}{10.00}{28}
\emline{90.00}{10.00}{29}{90.00}{20.00}{30}
\end{picture}

\noindent Then discrete concavity is equivalent to validity of
three types of ``rhombus'' inequalities. The inequalities require
that sum at two points of drawn diagonal is greater or equal to
the sum at two points of non-drawn diagonal.

 (i)   $f(i,j)+f(i+1,j+1)\ge f(i+1,j)+f(i,j+1)$;

(ii)  $f(i,j+1)+f(i+1,j+1)\ge f(i+1,j+2)+f(i,j)$.

(iii)  $f(i+1,j)+f(i+1,j+1)\ge f(i,j)+f(i+2,j+1)$;

Note, that if only the requirement (i) is valid, then a function
is called {\em supermodular}. If a function is supermodular and
the requirement (ii) is valid, then the function is discrete
concave on every vertical strip of the unit length, and we call
such functions  {\em vertically-strip concave} ($VS$-concave).
Analogously, if (i) and (iii) are valid, a function is called {\em
horizontally strip-concave}  ($HS$-concave).

Mostly, we will be interested in functions on the triangle grid
with the vertexes  $(0,0)$, $(0,n)$, $(n,n)$; denoted by $\Delta
_n$. On the next picture we depicted the grid $\Delta _4$.

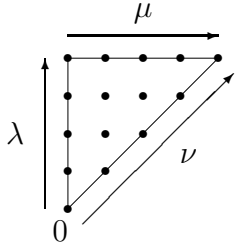
\unitlength=1.00mm \special{em:linewidth 0.4pt}
\linethickness{0.4pt}
\begin{picture}(100.00,35.00)(-10,0)
\emline{50.00}{5.00}{1}{50.00}{25.00}{2}
\emline{50.00}{25.00}{3}{70.00}{25.00}{4}
\emline{70.00}{25.00}{5}{50.00}{5.00}{6}
\put(47.00,5.00){\vector(0,1){20.00}}
\put(50.00,28.00){\vector(1,0){20.00}}
\put(52.00,3.00){\vector(1,1){20.00}}
\put(49.00,2.00){\makebox(0,0)[cc]{$0$}}
\put(43.00,15.00){\makebox(0,0)[cc]{$\lambda$}}
\put(60.00,31.00){\makebox(0,0)[cc]{$\mu$}}
\put(66.00,12.00){\makebox(0,0)[cc]{$\nu$}}
\put(50.00,5.00){\circle*{1.00}} \put(50.00,10.00){\circle*{1.00}}
\put(50.00,15.00){\circle*{1.00}} \put(50.00,20.00){\circle*{1.00}}
\put(50.00,25.00){\circle*{1.00}} \put(55.00,25.00){\circle*{1.00}}
\put(60.00,25.00){\circle*{1.00}} \put(65.00,25.00){\circle*{1.00}}
\put(70.00,25.00){\circle*{1.00}} \put(55.00,20.00){\circle*{1.00}}
\put(60.00,20.00){\circle*{1.00}} \put(65.00,20.00){\circle*{1.00}}
\put(55.00,15.00){\circle*{1.00}} \put(60.00,15.00){\circle*{1.00}}
\put(55.00,10.00){\circle*{1.00}}
\end{picture}

Consider a discrete concave function $f$ on the grid
$\Delta_n$\footnote{Such a discrete concave function was called a
{\em hive} in \cite{KT,oct}.} and consider its restriction to each
side of the triangle:  the left-hand side, the top of the triangle
and the hypotenuse. Specifically, we orient these sides as
depicted on the previous picture and consider increments of the
function on each unit segment. Then, increments along the
left-hand side constitute an $n$-tuple
$$
\lambda (1)=f(0,1)-f(0,0), \lambda (2)=f(0,2)-f(0,1),...,\lambda
(n)=f(0,n)-f(0,n-1).
$$
It is easy follows from the rhombus inequalities of the type (i)
and (iii) that
$$
\lambda (1)\ge \lambda (2)\ge ...\ge \lambda (n).
$$
Analogously, we define  $n$-tuple $\mu $ ($\mu
(i)=f(i,n)-f(i-1,n)$, $i=1,...,n$) and $\nu $ ($\nu
(k)=f(k,k)-f(k-1,k-1)$, $k=1,...,n$), which are also decreasing
tuples. We call these $n$-tuples {\em increments} of the function
$f$ on the corresponding sides of the triangle grid. Obviously,
the increments are invariant under adding a constant to $f$.
Therefore, we have to consider functions modulo adding a constant
or to require $f(0,0)=0$.

Let us briefly say about main roles of discrete concave functions
in combinatorics and representation theory. We let to denote
$DC_n(\lambda ,\mu ,\nu )$ the set of discrete concave functions
on the grid  $\Delta _n$ with increments  $\lambda$, $\mu $, $\nu
$. This set is a polytope (probably empty) in the space of all
functions on $\Delta _n$. If this polytope is non-empty, when the
$n$-tuples $\lambda ,\mu ,\nu $ are decreasing and there holds
$|\lambda |+|\mu |=|\nu |$. For  $n>2$, we need more relations in
order to get a non-empty $DC_n(\lambda ,\mu ,\nu )$. The necessary
and sufficient conditions for non-emptyness of  $DC_n(\lambda ,\mu
,\nu )$ (so-called Horn inequalities)  are in  т \cite{KT}, see
also \cite{Fhorn,Kar,stekl}. Moreover, $DC_n(\lambda ,\mu ,\nu )$
is non-empty if and only if there exist Hermitian matrices  $A$
and $B$, such that  $A$, $B$, $A+B$ have spectra $\lambda ,\mu
,\nu $, respectively  (a solution to the Horn problem).

We let to denote $DC_n^{\mathbb Z}(\lambda ,\mu ,\nu )$ the set of
integer-valued discrete concave functions on the grid  $\Delta
_n$, of course the tuples $\lambda$, $\mu $, $\nu $ have to be
integer-valued as well. The cardinality of this set coincides with
the Littlewood-Richardson coefficient, the multiplicity of the
irreducible representation $V_{\nu}$ (of $GL(n)$) in the tensor
product irreps $V_{\lambda}\otimes V_{\mu}$. In Section 7 we will
be more specific on this issue.\medskip

\section{Functions on {\bf 3D}-grids}

For the purpose of this section, it is convenient to consider the
octahedron recurrence with the propagation vector $(-1,1,1)$ and
locate the initial data at the qudrants  $OXZ$ ({\em ground}) and
$OXY$ ({\em front wall}). The modular flats take the form $x=a$,
$y=b$, $z=c$ and $x+y+z=d$, where  $a, b, c, d\in \mathbb Z $. If
we cut $\mathbb R^3$ by these planes, we get a decomposition of
$\mathbb R^3$ into primitive tetrahedrons and octahedrons. All
octahedrons are parallel, that is one can be obtained by an
integer translation of another. Each octahedron has three
diagonals parallel to vectors $(1,1,-1)$, $(1,-1,1)$ and
$(-1,1,1)$, respectively, and corresponding three pairs of
antipodal vertexes.

The diagonal being parallel to the propagation vector $(-1,1,1)$,
we call the {\em mail} diagonal. The OR leads us to the following
notion.
\medskip

{\bf Definition.} A function $F:\mathbb Z^3\to\mathbb
R\cup\{-\infty\}$ is said to be {\em polarized}, if, for any
primitive octahedron, sum of values of  $F$ at the vertexes of the
main diagonal is equal to the maximum of the sum of values of  $F$
at the antipodal vertexes of two others diagonals.\medskip

We denote $\Delta_n  (OXYZ)$ the three-dimensional grid,
constituted from the non-negative integer points   $(x,y,z)$, such
that $x+y+z\le n$. It is easy to see that, for any initial data
given at the {\em ground } $\Delta_n (OXY)$ and the {\em front
wall} $\Delta_n (OXZ)$, there exists a unique polarized function
with domain  $\Delta_n (OXYZ)$ and these given values. This is
done by the OR. However, we can set initial data at the  {\em
shadow wall} $\Delta _n (OYZ)$ and the {\em slope wall} $\Delta _n
(XYZ)$ and get a polarized function. In that case, we have to
apply the OR with the reverse propagation vector $(1,-1,-1)$. On
the next Picture we depicted the  tetrahedron co$\Delta_n  (OXYZ)$
with the direction of the OR propagation.

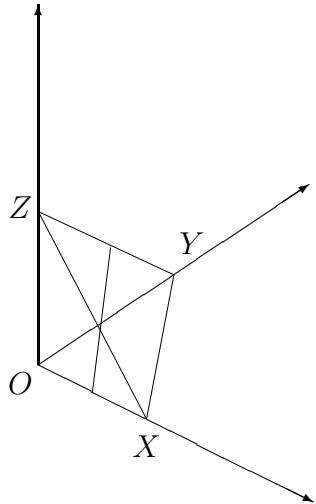
\unitlength=0.60mm \special{em:linewidth 0.4pt}
\linethickness{0.4pt}
\begin{picture}(82.00,141.00)
\put(21.00,61.00){\vector(3,2){60.00}}
\put(21.00,61.00){\vector(0,1){80.00}}
\put(21.00,61.00){\vector(2,-1){61.00}}
\emline{51.00}{81.00}{1}{21.00}{95.00}{2}
\put(17.00,57.00){\makebox(0,0)[cc]{$O$}}
\put(17.00,96.00){\makebox(0,0)[cc]{$Z$}}
\put(55.00,88.00){\makebox(0,0)[cc]{$Y$}}
\put(45.00,43.00){\makebox(0,0)[cc]{$X$}}
\emline{51.00}{81.00}{3}{45.00}{49.00}{4}
\emline{45.00}{49.00}{5}{21.00}{95.00}{6}
\emline{33.00}{55.00}{7}{37.00}{87.00}{8}
\end{picture}

The fundamental property of the octahedron recurrence is that if
the initial data (at the ground and the front wall) are discrete
concave function, then the corresponding polarized function on the
grid $\Delta_n (OXYZ)$ is a kind of three-dimensional discrete
concave function. Without going in details of discrete concave
functions in $\mathbb Z^n$, we give notions appropriate for this
paper.

Discrete concavity on {\bf 2D}-grids is equivalent to fulfil three
kinds of rhombus inequalities. In dimension 3, we have four kinds
of modular flats. In each such a 2-dimensional flat we have
rhombuses, which corresponds to triangular decomposition of the
flat by cutting it by three others kinds of modular flats. We have
to require validity of rhombus inequality for each such a rhombus:
the sum of values at the ``short'' diagonal is greater or equal to
the sum at the ``long'' diagonal.\medskip

{\bf Definition}. A function $F:\mathbb Z^3\to\mathbb
R\cup\{-\infty\}$ is a {\em polarized discrete concave} function
if $F$ is polarized and all kinds of rhombus inequalities in each
modular flat are fulfilled. \medskip

Let us denote by $PDC_n$ the set of polarized discrete concave
functions on the three-dimensional grid $\Delta _n
(OXYZ)$.\medskip

{\bf Theorem 1}.  {\em  Let $F$ be a polarized function on the
three-dimensional grid $\Delta _n (OXYZ)$. Suppose the restriction
of $F$ to the ground face $\Delta_n (OXY)$ and to the { front
wall} face $\Delta_n (OXZ)$ are 2-dimensional discrete concave
functions. Then  $F\in PDC_n$. }\medskip

For a proof see \cite{ass}. Note, that this theorem is equivalent
to the following corollary (a sketch of proof of which is also in
\cite{h-k}). \medskip

{\bf  Corollary 1}. {\em  If the restrictions of a polarized
function to the ground and the front wall faces are discrete
concave, then the restriction to the shadow wall and the slope
wall are also discrete concave.}\medskip

{\em Proof}. In fact, any rhombus located on the slope or shadow
wall is also a rhombus for three-dimensional grid, and therefore,
the corresponding rhombus inequality is valid. $\Box$\medskip

{\bf Corollary 2} {\em Let the restriction of a polarized function
to the ground be  discrete concave and the restriction to the
front wall be $HS$-concave. Then the restrictions to two other
faces are $HS$-concave.}
\medskip

{\em Proof}. In fact, we can add to  $F$ an appropriate function
$\varphi (z)$ of the vertical variable  $z$, in order to get a
discrete concave function on the front wall. $F+\varphi (z)$ is
not changed on the ground, therefore by Corollary 1,  $F+\varphi
(z)$ is discrete concave at the other two wall, therefore $F$ is
$HS$-concave on these walls. $\Box$\medskip

{\bf Corollary 3}. {\em  Suppose the restriction to the ground of
a polarized function is $HS$-concave and the restriction to the
front wall is $VS$-concave (here we consider horizontal being
parallel to the segment $XY$). Then $F$ is $VS$-concave on the
shadow wall. }\medskip

{\em Proof}. As above, having add to $F$ an appropriate separable
function on variables $x$ and $y$, we get a polarized function
$G=F+\varphi (x)+\psi (y)$, which will be discrete concave on the
ground and the front wall. By Corollary 1, $G$ is discrete concave
on the shadow wall. Therefore, $F$ is $VS$-concave on this wall.
$\Box$\medskip

{\bf Corollary 4}. {\em Suppose a polarized function $F$ is
discrete concave on the ground and $VS$-concave on the front wall.
Then $F$ is discrete concave on the shadow wall.}\medskip

{\em Proof}. In fact, having add an appropriate function on $x$ to
$F$, we get a discrete concave function on the front wall. On the
ground this function will be also discrete concave. But this
function remains the same on the shadow wall, and by Corollary 1
the function on this wall is discrete concave. $\Box$\medskip

Now let us consider the polarized functions (or the octahedron
recurrence) on the prism  $\Delta_n (OXY)\times \{0,1,...,m\}$
(see next Picture).

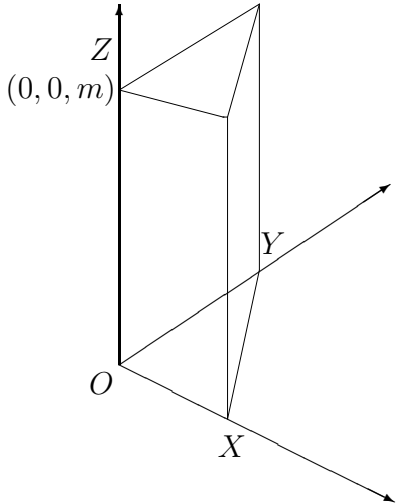
\unitlength=0.60mm \special{em:linewidth 0.4pt}
\linethickness{0.4pt}
\begin{picture}(82.00,141.00)
\put(21.00,61.00){\vector(3,2){60.00}}
\put(21.00,61.00){\vector(0,1){80.00}}
\put(21.00,61.00){\vector(2,-1){61.00}}
\put(17.00,57.00){\makebox(0,0)[cc]{$O$}}
\put(17.00,131.00){\makebox(0,0)[cc]{$Z$}}
\put(55.00,88.00){\makebox(0,0)[cc]{$Y$}}
\put(46.00,43.00){\makebox(0,0)[cc]{$X$}}
\emline{45.00}{49.00}{1}{52.00}{82.00}{2}
\emline{52.00}{82.00}{3}{52.00}{141.00}{4}
\emline{52.00}{141.00}{5}{21.00}{122.00}{6}
\emline{21.00}{122.00}{7}{44.00}{116.00}{8}
\emline{44.00}{116.00}{9}{44.00}{116.00}{10}
\emline{44.00}{116.00}{11}{45.00}{116.00}{12}
\emline{45.00}{116.00}{13}{52.00}{141.00}{14}
\emline{45.00}{116.00}{15}{45.00}{49.00}{16}
\put(8.00,122.00){\makebox(0,0)[cc]{$(0,0,m)$}}
\end{picture}

We have to set functions equals $-\infty$ at points outside the
prism. Therefore, on the non-modular face $\Delta_n (XY)\times
\{0,1,...,m\}$, a polarized function  $F$ has to be a separable
function (on variables $x+y$ and $z$). In other words, for any
``primitive'' quadrat on this wall, the sum of values at the
opposite pairs of vertexes coincide.\medskip

{\bf  Corollary 5} {\em  Let $F$ be a polarized function on the
prism  $\Delta_n (OXY)\times \{0,1,...,m\}$. Suppose the
restriction of $F$ to the ground face $\Delta_n (OXY)\times \{0\}$
and the restriction to the front wall  $\Delta _n (OX)\times
\{0,1,...,m\}$ are discrete concave functions. Then $F$ is
polarized discrete concave function on the prism (and, in
particular, $F$ is discrete concave on the shadow wall $\Delta _n
(OY)\times \{0,1,...,m\}$ and on the ceiling  $\Delta_n
(OXY)\times \{m\}$). }\medskip

{\em Proof}. It is easy to see that it suffices to prove the
corollary in the case $m=2$.

In the beginning we consider the case $m=1$. Let us extend the
ground to the size of $n+1$, that is we add to $\Delta_n(OXY)$ new
points $(n+1,0,0),...,(0,n+1,0)$. Let us  extend  $F$ to these
points such that we  get a discrete concave function on the
extended ground $\Delta_{n+1} (OXY)\times \{0\}$ and a discrete
concave function on the ``extended'' front wall. We can always do
that by setting small values ($<<0$) to these points. Let us
denote $\tilde F$ such an extension. By Theorem 1, the function
$\tilde F$ is a polarized discrete concave function. We claim,
that the restriction of this function to the prism is a polarized
discrete concave function. In fact, it suffices to check that
$\tilde F$ coincides with  $F$ on the non-modular face $\Delta_n
(XY)\times \{0,1\}$. But this holds since we assigned small values
to the new points. Thus  $F$ and $\tilde F$ coincide on the prism.
Since  $\tilde F$ is discrete concave function,  $F$ is discrete
concave too.

Let us move to the case $m=2$. We have to check all rhombus
inequalities for all rhombuses in the prism $\Delta _n (OXY)\times
\{0,2\}$. Let us first consider the rhombuses of the vertical size
2. It is easy to see that these rhombuses belong to the
tetrahedron of size $n+1$. Then the corresponding rhombus
inequality is valid, since they are valid for  $\tilde F$. Other
rhombuses  are located either in the prism  $\Delta_n (OXY)\times
\{0,1\}$, or in the prism  $\Delta _n (OXY)\times \{1,2\}$. For
the first prism, the corresponding inequality follows due to the
above case with $m=1$. Moreover, we get that $F$ is discrete
concave on the triangle  $\Delta_n (OXY)\times \{1\}$. Now, again
applying the case $m=1$ to the prism $\Delta _n (OXY)\times
\{1,2\}$, we get validity of rhombus inequalities in this prism.
$\Box$\medskip

Using similar reasonings one can get the following\medskip

{\bf  Corollary 6} {\em Suppose a polarized function  $F$ on the
prism  $\Delta _n (OXY)\times \{0,1,...,m\}$ has discrete concave
restrictions to the ceiling $\Delta _n (OXY)\times \{m\}$ and the
shadow wall $\Delta _n (OY)\times \{0,1,...,m\}$. Then  $F$ is a
polarized discrete concave function and its restrictions to the
front wall $\Delta _n (OX)\times \{0,1,...,m'\}$ and the ground
$\Delta_n (OXY)\times \{0\}$ are discrete concave functions. }

\section{Arrays}

In this section, we introduce another key player of a game --
arrays. Consider a rectangle  $[0,n]\times [0,m]$ on the plane
with natural  $n$ and $m$, constituted from unit squares with the
centers at the points $(i-1/2,j-1/2)$,  $i=1,...,n$, $j=1,...,m$,
we call such squares boxes. An array is a filling of each box
$(i,j)$ with a non-negative ``mass'' $a(i,j)$.

To each array $a$ we associate a function $f=f_a$ on the
rectangular grid  $\{0,1,...,n\}\times \{0,1,...,m\}$ by setting
to the point  $(i,j)$ the value
$$
f_a(i,j)=\sum_{i'\le i,\,j'\le j} a(i',j').
$$
In other words, this value is equal to the mass of all boxes to
the south-west from the point $(i,j)$. This is a reason to denote
by  $\int \! \! \int a$ the function $f_a$. On the bottom and the
left boundary of the rectangle the function equals $0$. For other
$(i,j)$, we obviously have
$$
f(i,j)-f(i-1,j)-f(i,j-1)+f(i-1,j-1)=a(i,j).
$$
From this  $a(i,j)$ might be understand as the mixed derivative of
$f$ ($a=\partial\partial f$), or as a break of  $f$ along the
common edge $[(i-1,j-1),(i,j)]$ of two affinity areas. Since
$a(i,j)\ge 0$, the function $\int \! \! \int a$ is supermodular.

Here we collect some notions and results on arrays which are of
use in the paper (for details see  \cite{umn}).\medskip

\begin{enumerate}
\item  For each  $j=1,...,m-1$, there is an operation  $D_j$ on
the set of arrays, which acts on a given array by moving down
(vertically) an amount of mass ( $\le 1$) from a box $(i,j)$ to
$(i,j-1)$(for a definition such a $i$ see \cite{umn} (see also
\cite{kyoto})). For any $a$, starting from some power
$\epsilon_j(a)$, there holds $D_j^{\epsilon
(a)+1}(a)=D_j^{\epsilon (a)}(a)$, we denote ${\bf
D}_j=D_j^{+\infty}$, that is ${\mathbf D}_j(a)=D_j^{\epsilon
(a)}(a)$.

\item  If, for an array $a$, there holds $D_j(a)=a$  (or ${\bf D}_ja=a$), then $a$
is called ${\bf D}_j$-tight. Equivalently, this means that the
function $f=\int \!\!\int a$ satisfies the inequalities
$$
              f(i-1,j)+f(i,j)-f(i-1,j-1)-f(i,j+1)\ge 0
$$
for all $i=1,...,n$. In other words, the rhombus inequalities of
type (ii) (see Section 3) hold true for the rhombuses which are
cut by the line  $y=j$.

\item If an array $a$ is  ${\bf D}_j$-tight for all
$j=1,...,m-1$, then $a$ is said to be  {\bf D}-{\em tight}. For a
{\bf D}-{\em tight} array $a$, the corresponding function
$f_a=\int\!\!\int a$ is a $VS$-concave function.

It is clear that  we can condense any array to a {\bf D}-tight
(for example, by applying $(\mathbf D_1\ldots\mathbf D_{m-1})^m$,
but this is only one of ways). Moreover, for each $a$, such a {\bf
D}-tight array is defined uniquely and we let to denote it by
${\bf D}a$. Since, descending massed due to the operations $D_j$,
does not change the vector of column sums (masses), the values of
the functions $\int \! \! \int {\bf D}a$ and $\int \! \! \int a$
coincide at the top boundary of the rectangle, i.e. at the points
with $y=m$. At the right boundary, i.e., for $x=n$, the values are
different (for non $\mathbf D$-tight arrays). The increments of
the function $\int \! \! \int {\bf D}a$ along the right side we
let to denote by $\lambda_1 ,...,\lambda_{m}$. Due to
$VS$-concavity of the function $\int \! \! \int {\bf D}a$, we have
the inequalities $\lambda _1 \ge \lambda _2 \ge ...\ge \lambda_{m}
\ge 0$. This m-tuple  $\lambda $ we call  ${\bf D}$-{\em shape} of
$a$.

Obviously, for a   $\mathbf{D}$-tight array $a$, the ${\bf
D}$-{\em shape} of $a$ coincides with the vector of its row sums.

Integer-valued $\mathbf D$-tight arrays are in a natural bijection
with semistandard Young tableaux. For details see \cite{umn}, and
here we explain this bijection by an example.

{\bf Example}. Consider the following $\mathbf D$-tight $4\times
3$ array
$$
\pmatrix{0&0& 0&3\cr 0&4&0&4\cr 5&1&2&4}
$$
To get the corresponding semi-standard Young tableaux we have to
read this array from left to right and from bottom to top. Reading
a row gives us filling of the corresponding row in the Young
tableau, the mass $a(i,j)$ exhibits the multiplicity of
repetitions of the letters $i$ in the $j$-th row of the Young
tableaux (we consider the French style of drawing Young diagrams
and tableaux, that is the Young diagram for a partition
$\lambda_1\ge \lambda_2\ge\ldots\ge \lambda_n\ge 0$ is a
collection of boxes in the grid  $\mathbb N\times \mathbb N$ with
north-east corners $(i,j)$ such that $j\le n$, $i\le \lambda_j$,
and a Young tableau is a filling of the diagram from some alphabet
increasingly along each row (from left to right) and strictly
increasing from bottom to top). Thus, for the above array, we get
$$
\begin{array}{ccccccccccccc}
4&4&4\\ 2&2&2&2&4&4&4&4\\
1&1&1&1&1&2&3&3&4&4&4
\end{array}
$$
\medskip

\item Using the transposition (with respect to the diagonal $a^T(i,j)=a(j,i)$)
we define the operations $L_i$ ${\bf L}_i$ (which translate masses
to the left along a row), $L_i(a)=(D_i(a^T))^T$. Using the
operations $L_i$, $i=1,\ldots, n-1$, we can  condense to the left
any array $a$, and get  {\bf L}-tight array  ${\bf L}a$. The
function $\int \!  \! \int {\bf L}a$ is  an $HS$-concave function.
In particular, the increments of this function along the top side
is a decreasing tuple, the column sum of the array ${\bf L}a$,
that is the  {\bf L}-shape of  $a$.

\item The key claim of theory of arrays is that the
operations $L_i$ and $D_j$ commute for any $i$, $j$ (\cite{umn},
Theorem  4.2). Here we present some consequences of this
commutation property
\medskip

a) For any array $a$,  the {\bf D}-shape of $a$ coincides with the
{\bf L}-shape of $a$, and, thus, this tuple is the {\em shape} of
$a$.\medskip

b) The bijection theorem (or modified  $RSK$ correspondence,
\cite{umn}, Theorem  6.2): suppose we are given a  {\bf D}-tight
array  $d$ and an {\bf L}-tight array  $l$, such that $d$ and $l$
have the same shape, then there exists a unique array  $a$, such
that $d={\bf D}a$ and $l={\bf L}a$ hold true.
\end{enumerate}
Our next task is to obtain this modified  $RSK$ bijection using
the octahedron recurrence.

\section{Functional form of $RSK$}

Let us turn back to Picture 2 and locate the zero function at the
face  $OEA$; at the slope rectangle $OABC$ we  locate the function
$f_a=\int \!\! \int a$. Specifically, we assign the value $\int
\!\!\int a(i,j)$) to the point  $(i,j,j)$. Now, we propagate these
data by the octahedron recurrence to the prism. From Corollary 5
(Section 4), we get a VS-concave function at the top face
rectangle $EABD$ and HS-concave function at the right face
triangle $CDB$ (the vertical is $y$-axe in the first case and
$z$-axe in the second case). Moreover, we get the function $\int
\! \! \int {\bf D}a$ {\em at the tope face and the function $\int
\! \! \int {\bf L}a$ at the right face of the prism}
(specifically, the restriction of this function to this triangle).

Namely, we state\medskip

{\bf Theorem 2}. {\em  Let $F$ denote a function on the prism
obtained by the octahedron recurrence from  the following initial
data: the zero values at the faces $OEDC$ and $OEA$, and $\int \!
\! \int a$ at $OABC$. Then $F(i,j,m)=(\int \! \! \int {\bf
D}a)(i,j)$ for all $i,j$, and $F(n,j,k)=(\int \!  \! \int {\bf
L}a)(j,k)$ for $0\le j\le k\le m$.}\medskip

Let us note, that the latter two functions coincide at the edge
$DB$. That is  $(\int \!  \! \int {\bf D}a)(n,j)=(\int \!  \! \int
{\bf L}a)(j,m)$ for all  $j$. In fact, this is $\int $ from the
shape of  $a$.\medskip

{\bf Remark}. Of course, we can consider a propagation in the
reverse direction. Specifically, assume we are given a function
$f$ on the top face  $EABD$ and a function  $g$ on the triangle
$CDB$. Suppose there hold

     a) $f$ is $VS$-concave and equals $0$ at the edges
$EA$ and $ED$;

     b) $g$ is $HS$-concave and equals  $0$ at the edge
$CD$;

     c) the functions $f$ and $g$ coincide at  the edge $DB$.

\noindent Then having apply the OR (with the propagation vector
$(-1,0,-1)$) for these data, we get a pair of functions on the
triangle $OEA$ and the slope rectangle $OABC$. Due to Corollary 6,
we get a discrete concave function on  $OEA$ and a supermodular
function on  $OABC$. Moreover, we get the identically zero
function on the triangle $OEA$. This is because this function
equals $0$ at the edge  $EA$ (from the item a)) and at the edge
$OE$ (this follows from  b) and separability of the OR on the
non-modular face). But these boundary values force nullity of the
discrete concave function.

Now, we get an array $a$ as the mixed derivatives of the
supermodular function on $OABC$. It is clear that  $f=\int\!\!\int
a$ and $g=\int\!\!\int a$. Thus, this octahedron recurrence
provides us with a functional form of the modified RSK (item 5 of
Section 5). An advantage of this form of RSK is that the direct
and inverse bijections are done symmetrically. \medskip

{\em  Proof of Theorem 2}. The main case of the proof is the case
of an array with two rows, that is  $m=2$. We denote the masses in
the bottom row by $a(i,1)$, and in the top row by $a(i,2)$,
$i=1,\ldots, n$. Therefore, for $i=1,\ldots, n$, we have
$f(i,0)=0$,
$$
                    f(i,1)=a(1,1)+...+a(i,1),
$$
$$
                 f(i,2)=f(i,1)+a(1,2)+...+a(i,2).
$$
This follows from the definition of $\int \!\!\int a$. Denote
$a'={\bf D}a$ and  $f'=\int\!\! \int a'$. Then $f'$ coincides with
$f$ for $j=0$ (both functions  $=0$) and for $j=2$. Differences
could occur only at points with $j=1$, because some mass moves
from the first level to the ground level. We get in \cite{umn},
formula $(***)$, the following  $ f'(i,1)=f(i,1)+\max(\beta _1
,...,\beta_i )$, where
$$\beta _i =a(1,2)+...+a(i,2)-a(1,1)-...-a(i-1,1)=$$
$$
f(i,2)-f(i,1)-f(i-1,1)+f(i-1,0).$$ Now, taking into account $
f(i,2)=f'(i,2)$, and
$$
\max(\beta _1 ,...,\beta _i )=\max(\max(\beta _1 ,...,\beta_{i-1}
),\beta _i ), $$ and
$$\max(\beta _1
,...,\beta_{i-1})=f'(i-1,1)-f(i-1,1), $$ we obtain
$$
f'(i,1)=
$$
$$f(i,1)+\max[f'(i-1,1)-f(i-1,1),f'(i,2)-f(i,1)-f(i-1,1)+f(i-1,0)]=
$$
$$
\max[f'(i-1,1)+f(i,1),f'(i,2)+f(i-1,0)]-f(i-1,1).
$$
That is the octahedron recurrence indeed.\medskip

We claim that in the prism $OEACDB$ at the height $z=k$ is located
the function  $\int \! \! \int {\bf D}(a_k )$, where the array
$a_k$ is obtained from  $a$ by omitting the rows with
$j=k+1,...,m$. In fact, suppose this claim is true for some $k$,
and let us check it for  $k+1$. So, we are given a  {\bf D}-tight
array $d_k ={\bf D}(a_k )$, and the row  $a(\cdot ,k+1)$ from the
array $a$. From the above case with $m=2$ follows that the
octahedron recurrence lifting, from the height $k$ to the height
$k+1$ in the prism, corresponds to the product of the condensation
operations  ${\bf D}_1 ...{\bf D}_k$. The formula 6.5 in
\cite{umn} demonstrates exactly this claim.\medskip

Thus at the top face, we get the function $\int \!  \! \int {\bf
D}(a)$.

Now, we have to show that at the right side face we get the
function $\int \! \! \int {\bf L}(a)$. In fact, from the above
claim, on the segment $\{(n, 0, k),\cdots , (n, k,k)\}$,  we have
the function  $\int \! \! \int d_k =\int \! \! \int {\bf D}(a_k
)$. That is the integral  $\int $ of the shape of the array $d_k$,
or, equivalently, the integral of the shape of the array  $a_k$
(see \cite{umn}, 5.10). But the shape of  $a_k$ equals the shape
of ${\bf L}(a_k )={\bf L}(a)_k$. For any {\bf L}-tight array, its
shape coincides with the column sum vector. This completes the
proof. $\Box$\medskip

Let us illustrate this theorem by an example. Consider the
following array $a=\begin{array}{ccc} 1&2&2\cr 1&1&5\cr 2&3&1
\end{array}$. The corresponding supermodular function
$f_a=\begin{array}{ccc} 4&10&18\cr 3&7&13\cr 2&5&6
\end{array}$
is located at the face  $OABC$; the values of the polarized
function $F$ are depicted on the next Picture\medskip

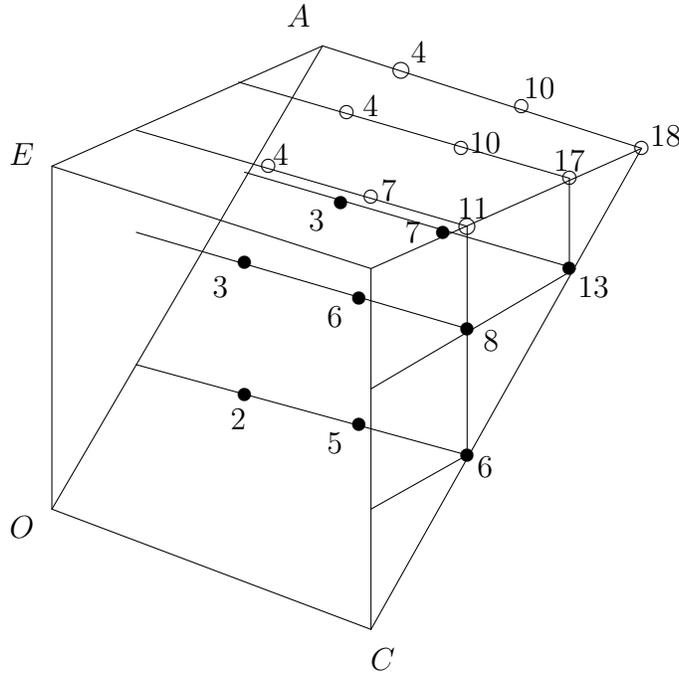
\unitlength=0.8mm \special{em:linewidth 0.4pt}
\linethickness{0.4pt}
\begin{picture}(119.00,122.00)
\emline{17.00}{40.00}{1}{17.00}{97.00}{2}
\emline{17.00}{97.00}{3}{70.00}{80.00}{4}
\emline{70.00}{80.00}{5}{70.00}{20.00}{6}
\emline{70.00}{20.00}{7}{17.00}{40.00}{8}
\emline{70.00}{80.00}{9}{115.00}{100.00}{10}
\emline{115.00}{100.00}{11}{62.00}{117.00}{12}
\emline{62.00}{117.00}{13}{17.00}{97.00}{14}
\emline{115.00}{100.00}{15}{70.00}{20.00}{16}
\emline{62.00}{117.00}{17}{17.00}{40.00}{18}
\emline{70.00}{60.00}{19}{104.00}{80.00}{20}
\emline{104.00}{80.00}{21}{49.00}{96.00}{22}
\emline{70.00}{40.00}{23}{86.00}{49.00}{24}
\emline{86.00}{49.00}{25}{31.00}{64.00}{26}
\emline{86.00}{49.00}{27}{86.00}{87.00}{28}
\emline{86.00}{87.00}{29}{31.00}{103.00}{30}
\emline{103.00}{80.00}{31}{103.00}{95.00}{32}
\emline{103.00}{95.00}{33}{48.00}{111.00}{34}
\emline{86.00}{70.00}{35}{31.00}{86.00}{36}
\put(86.00,49.00){\circle*{2.00}}
\put(68.00,54.00){\circle*{2.00}}
\put(49.00,59.00){\circle*{2.00}}
\put(86.00,70.00){\circle*{2.00}}
\put(68.00,75.00){\circle*{2.00}}
\put(49.00,81.00){\circle*{2.00}}
\put(103.00,80.00){\circle*{2.00}}
\put(82.00,86.00){\circle*{2.00}}
\put(65.00,91.00){\circle*{2.00}} \put(86.00,87.00){\circle{2.83}}
\put(70.00,92.00){\circle{2.00}} \put(53.00,97.00){\circle{2.00}}
\put(115.00,100.00){\circle{2.00}}
\put(95.00,107.00){\circle{2.00}}
\put(75.00,113.00){\circle{2.83}}
\put(103.00,95.00){\circle{2.00}}
\put(85.00,100.00){\circle{2.00}}
\put(66.00,106.00){\circle{2.00}}
\put(78.00,116.00){\makebox(0,0)[cc]{$4$}}
\put(98.00,110.00){\makebox(0,0)[cc]{$10$}}
\put(119.00,102.00){\makebox(0,0)[cc]{$18$}}
\put(70.00,107.00){\makebox(0,0)[cc]{$4$}}
\put(89.00,101.00){\makebox(0,0)[cc]{$10$}}
\put(103.00,98.00){\makebox(0,0)[cc]{$17$}}
\put(55.00,99.00){\makebox(0,0)[cc]{$4$}}
\put(73.00,93.00){\makebox(0,0)[cc]{$7$}}
\put(87.00,90.00){\makebox(0,0)[cc]{$11$}}
\put(48.00,55.00){\makebox(0,0)[cc]{$2$}}
\put(64.00,51.00){\makebox(0,0)[cc]{$5$}}
\put(89.00,47.00){\makebox(0,0)[cc]{$6$}}
\put(107.00,77.00){\makebox(0,0)[cc]{$13$}}
\put(90.00,68.00){\makebox(0,0)[cc]{$8$}}
\put(64.00,72.00){\makebox(0,0)[cc]{$6$}}
\put(45.00,77.00){\makebox(0,0)[cc]{$3$}}
\put(61.00,88.00){\makebox(0,0)[cc]{$3$}}
\put(77.00,86.00){\makebox(0,0)[cc]{$7$}}
\put(12.00,37.00){\makebox(0,0)[cc]{$O$}}
\put(12.00,99.00){\makebox(0,0)[cc]{$E$}}
\put(58.00,122.00){\makebox(0,0)[cc]{$A$}}
\put(72.00,15.00){\makebox(0,0)[cc]{$C$}}
\end{picture}

The values of the function $F$ at the integer points of top face
$EABD$ are $\begin{array}{ccc} 4&10&18\cr 4&10&17\cr 4&7&11
\end{array}$ and the corresponding SSYT is
$$
\begin{array}{cccccccccccc}
3& & & & & & & & & & &\cr 2&2&2&3&3&3& & & & & &\cr
1&1&1&1&2&2&2&3&3&3&3&
\end{array}.
$$
The values of $F$ at the face $CBD$ are $\begin{array}{ccc}
\phantom{6}& \phantom{13}&18\cr \phantom{6}&13&17\cr 6&8&11
\end{array}$
and the corresponding SSYT is
$$
\begin{array}{cccccccccccc}
3& & & & & & & & & & &\cr 2&2&2&2&2&3& & & & & &\cr
1&1&1&1&1&1&2&2&3&3&3&
\end{array}.
$$\medskip

Now we give consequences of this theorem for the associativity and
commutativity bijections.

\section{ Associativity bijection}

Let us briefly  recall the matter. Pick an integer  $n\ge 1$, let
to denote by  small Greek letters $\lambda$, $\mu $, $\nu $ etc.
partitions with  $n$-parts (that is an n-tuple $\lambda =(\lambda
_1 ,...,\lambda _n )$ of integers such that $\lambda _1 \ge
\lambda _2 \ge ...\ge \lambda _n\ge 0 $). To each partition
$\lambda $ is associated the Schur function $s_{\lambda }$ (see
\cite{M,umn}). While $\lambda$ runs over the set of all
$n$-partitions, the Schur functions constitute an additive basis
of the ring of symmetric functions on $n$ variables. Therefore,
the product  $s_{\lambda }s_{\mu }$ of the Schur functions can be
presented of the form
$$
s_{\lambda} s_{\mu}=\sum_{\nu} c_{\lambda,\mu }^{\nu } s_{\nu}.
                  $$
The structure constants $c_{\lambda ,\mu }^{\nu }$ are called the
{\em Littlewood-Richardson coefficients}. Since the ring of
symmetric functions is commutative and associative, the
Littlewood-Richardson coefficients satisfy the commutativity
$$
c_{\lambda ,\mu }^{\nu }=c_{\mu, \lambda }^{\nu }
$$
and associativity
$$
 \sum _{\sigma } c_{\lambda ,\mu }^{\sigma}  c_{\sigma ,\nu }^{\pi}
= \sum _{\tau}  c_{\mu ,\nu }^{\tau } c_{\lambda ,\tau }^{\pi }.
$$

Combinatorics learns us to seek for numbers underlying finite
sets, and for equalities look for natural bijections between the
corresponding sets. Littlewood and Richardson were the first who
give a combinatorial rule for computing LR-coefficients as the
cardinality of a set of special semistandard skew Young tableaux
(for example, from this, follows that these coefficients are
non-negative). Recently (\cite{BZ,Buch,stekl,oct}), these special
semistandard Young tableaux have been identified with
integer-valued discrete concave functions. Here we will give one
more interpretation of $c_{\lambda,\mu }^{\nu }$ as the
cardinality of the set of standard pairs of arrays $SP_{\mathbb
Z}(\lambda ,\mu ,\nu )$. In the array language, the commutativity
means equal cardinality of the sets $SP_{\mathbb Z}(\lambda ,\mu
,\nu )$ and  $SP_{\mathbb Z}(\mu ,\lambda ,\nu )$. Moreover, we
construct a natural bijection between these sets. Regarding the
associativity, we will provide a natural bijection between the
following sets
$$
\coprod _{\sigma }(SP_{\mathbb Z}(\lambda ,\mu ,\sigma )\times
SP_{\mathbb Z}(\sigma ,\nu ,\pi ))
$$
and
$$
\coprod _{\tau }(SP_{\mathbb Z}(\mu ,\nu ,\tau )\times SP_{\mathbb
Z}(\lambda ,\tau ,\pi )).
$$
The associativity bijection in \cite{oct} was constructed in terms
of hives (discrete concave functions) and the octahedron
recurrence played the main role. However, despite on an elegance
of the construction, a reason why it is natural (and even why it
is a bijection) was obscured. Recall that this bijection, which we
will call the {\em functional associativity bijection} provides a
bijection between sets
$$
\coprod _{\sigma }(DC^{\mathbb Z}_n(\lambda ,\mu ,\sigma )\times
DC^{\mathbb Z}_n(\sigma ,\nu ,\pi )) \mbox{ and } \coprod _{\tau
}(DC^{\mathbb Z}_n(\mu ,\nu ,\tau )\times DC^{\mathbb Z}_n(\lambda
,\tau ,\pi )),
$$
and it takes the following form. We have to locate a pair of
discrete concave functions $(f,g)$ from the first set at the faces
$OEA$ and $OAB$, respectively, of the tetrahedron $OEAB$ (see
Picture 2. $m=n$). Then, due to the OR, we get a pair of discrete
concave functions $(p,q)$ on the other two faces from the second
set (due to Theorem 1). The mapping $(f,g)\to (p,q)$ provides the
functional associativity bijection.\medskip

Using a relation between arrays and functions, and a natural
associativity bijection in terms of arrays, we provide a
justification for this construction. Namely, we prove  that the
associativity bijection in the arrays terms coincides with this
functional bijection. Analogously, for the case of the
commutativity bijection, we prove that the commutativity bijection
for arrays coincides with commutativity bijection in \cite{h-k}
and with two fundamental symmetries due to Pak and Vallexo
(\cite{PV}, Conjecture 1).\medskip

Here we consider square arrays of size  $n\times n$. Let us pick
an array  $a$, and consider the collection of arrays of the form
$Ta$, where $T$ is an arbitrary word in non-commutative variables
$D_j$ and  $U_j$, $j=1,\ldots, n-1$. This set constitute the {\em
orbit}  $O(a)$ of this array under action of the semi-group,
spanned by $D_j$ and  $U_j$, $j=1,\ldots, n-1$. In particular,
${\bf D}a \in O(a)$. It is not difficult to check (see
\cite{umn}), that each orbit contains a unique  {\bf D}-tight
array among its elements. Thus, to set an orbit is equivalent to
pick a  {\bf D}-tight array.

Furthermore,  if we take two  {\bf D}-tight arrays  $d$ and
$d'=R_id$ with some $i$ (or a pair of $\mathbf D$-tight arrays
which belong to an orbit under the action of the operations $R_i$
and $L_i$, $i=1,\ldots, n-1$), then the orbits, corresponding to
these arrays are isomorphic (in some sense). Thus, if we are
interested in orbits by modulo isomorphisms, we can consider the
orbit with simultaneously  {\bf D}-tight and $\mathbf L$-tight
array $d$ (or {\bf D}-tight and $\mathbf R$-tight, which is useful
sometimes). Such a bi-tight array is a partition $\lambda $
indeed. Specifically, for an $n$-tuple partition
$\lambda=(\lambda_1\ge ...)$, we denote by diag$(\lambda)$ an
array which contains the mass $\lambda_i$ at the diagonal box
$(i,i)$, $i=1,\ldots, n$, and zero masses at all other boxes. This
array  diag$(\lambda)$ is a bi-tight array, and any bi-tight array
takes such a form. We let to call such an orbit  $O(\lambda
)=O(diag(\lambda))$ the  {\em standard orbit} of shape $\lambda
$.\medskip

(Let us note, that irreducible (and polynomial) representations of
$GL(n)$ are also indexed by $n$-tuple partitions. Moreover, the
dimension of such a representation $V_{\lambda }$ coincides with
the cardinality of the orbit $O(\lambda )$. Therefore, one can
imagine the orbit $O(\lambda )$ as a skeleton of the irreducible
representation. Decomposition of an invariant set of arrays (under
the $DU$-action) into orbits corresponds to decomposition of a
representation into irreducibles. Moreover such a decomposition of
any invariant set is multiplicity-free.)\medskip

Let us consider two standard orbits  $O(\lambda )$ and $O(\mu )$.
We can form their tensor product  $O(\lambda )\otimes O(\mu )$. As
a set it is constituted of concatenations of arrays  $a\otimes b$,
$a\in O(\lambda)$, $b\in O(\mu)$, (that is an array of size  $2n
\times n$; the first $n$ columns come from $a$, and than come
columns of $b$). Since $\lambda $ and $\mu $ are {\bf L}-tight
arrays, the arrays $a$ and  $b$ are also {\bf L}-tight (commuting
$L$-operations and $U$-operations). The decomposition of
$O(\lambda )\otimes O(\mu )$ into orbits consists in
distinguishing $\mathbf D$-tight arrays of the form $a\otimes b$.
Thus, we come to the following
\medskip

{\bf Definition }. A {\em standard pair} is a pair of arrays
$(a,b)$ such that there holds  1) the arrays $a$ and $b$ are {\bf
L}-tight, and 2) the array  $a\otimes b$ is {\bf D}-tight. The
shape of  $a$ is the {\em starting } shape of the pair, the shape
of $b$ is the {\em intermediate} shape of the pair, and the shape
of the array $a\otimes b$, that is the vector of its row sums, is
the {\em final} shape.

We  let to denote  $SP(\lambda ,\mu ,\nu )$ the set of standard
pair with the starting shape $\lambda$, intermediate shape $\mu$
and the final shape $\nu$. The subset of integer valued standard
pair we mark with  ${\mathbb Z}$.
\medskip

Thus, the set of orbits in  $O(\lambda )\otimes O(\mu )$, which
are isomorphic to  $O(\nu )$, is identified to the set
$SP_{\mathbb{Z}}(\lambda ,\mu ,\nu )$. In \cite{umn}, 12.4, we
established that  $c_{\lambda,\mu}^{\nu}$ is equal to the
cardinality of the finite set  $SP_{\mathbb{Z}}(\lambda,\mu,\nu)$.
\medskip

Analogously, we can decompose the triple product $O(\lambda
)\otimes O(\mu )\otimes O(\nu )$ into the disjoint union of
orbits. An orbit of this decomposition is identified to a standard
triple  $(a,b,c)$ of arrays, such that  $a$, $b$ and $c$ are {\bf
L}-tight arrays (of shapes  $\lambda $, $\mu $ and $\nu $,
respectively), and the concatenated array  $a\otimes b\otimes c$
is {\bf D}-tight.

To a  standard triple we can correspond a couple of standard pairs
in two different ways, respectively to parenthesizes  $O(\lambda
)\otimes O(\mu )\otimes O(\nu )$. The parenthesizes $(O(\lambda
)\otimes O(\mu ))\otimes O(\nu )$, provides us with the pair
$(a,b)$ and the pair $({\bf L}(a\otimes b),c)$. Let us note, that
{\em the final shape of the first pair coincides with the starting
shape of the second pair}. Another parenthesizes $O(\lambda
)\otimes (O(\mu )\otimes O(\nu ))$ yields the pair $(a,{\bf
L}(b\otimes c))$ and the pair  ${\bf D}(b,c)$ (the latter is a
$(b',c')$, such that there holds $b'\otimes c'={\bf D}(b\otimes
c)$). It  is not difficult to check (see, for example \cite{umn})
that these pairs are standard indeed. (For example, because $b$
and $c$ are {\bf L}-tight, and since ${\bf D}$ and ${\bf L}_i$
commute, we get that $b'$ and  $c'$ are  {\bf L}-tight.) Let us
also note that, in the latter case, {\em the intermediate shape of
the first pair coincides with the final shape of the second pair}.
We call such couples of standard pairs {\em compatible couples}.
From the bijection theorem (see above the item 5b) follows that
such a compatible couple determines a standard triple $(a,b,c)$.
In fact, assume a compatible couple of standard pairs $(a,l)$ and
$(b',c')$ is given, and, in particular, there holds ${\bf D}l={\bf
L}(b',c')$. According to the bijection theorem there exists a pair
of arrays $b$ and  $c$, such that $l={\bf L}(b\otimes c)$ and
$b'\otimes c'={\bf D}(b\otimes c)$. Thus, we get a natural
associativity bijection for arrays
$$
\coprod_{\sigma}(SP(\lambda,\mu,\sigma)\times SP(\sigma,\nu,\pi)) \
\widetilde{\to} \ \coprod_{\tau}(SP(\mu,\nu,\tau)\times
SP(\lambda,\tau,\pi)),
$$
$$
((a,b),(\mathbf{L}(a\otimes b),c)) \mapsto ((a,\mathbf{L}(b\otimes
c)),\mathbf{D}(b,c)).
$$
Let us note, that in this construction the integer-validness of
$\lambda$, $\mu $, $\nu $, .... and the arrays do not play any
role, all is correct for arbitrary arrays.

\section{Comparing of two associativity bijections}

In the beginning, we present a natural bijection between the sets
$SP(\lambda ,\mu ,\nu )$ and  $DC_n(\lambda ,\mu ,\nu )$, where
$n$-tuples $\lambda,\mu,\nu$ are partitions. Namely, to a pair of
array  $(a,b)$ we assign the restriction of the function
$\int\!\!\int (a\otimes b)$ to the grid $n\le i\le n+j$, $0\le
j\le n$.\medskip

{\bf Proposition 2}. {\em  The above defined assignment provides a
bijection between sets $SP(\lambda ,\mu ,\nu )$ and $DC(\lambda
,\mu ,\nu )$. This assignment sends integer arrays into
integer-valued functions and vice versa.}
\medskip

{\em Proof}. Let $(a,b)$ be a standard pair. Since the array
$a\otimes b$ is {\bf D}-tight, the function $\int\!\!\int
(a\otimes b)$ is $VS$-concave on the rectangle  $2n\times n$.
Since the array $b$ is {\bf L}-tight, this function is
$HS$-concave on the square  $[n,2n]\times [0,n]$. Therefore, on
this square, the function $\int\!\!\int a\otimes b$ is discrete
concave. Furthermore, the increments of this function on the left
side of the square coincide with the vector of row sums of the
array $a$. Since, $a$ is a  {\bf D}-tight array, these increments
equal to the shape of  $a$, that is  $\lambda $. The increments on
the top side equal the column sums of  $b$. Since $b$ is a {\bf
L}-tight array, these sums equal the shape of $b$, that is  $\mu
$. Finally, increments on the right side  (``hypotenuse'')
coincide with the row sums of the array  $a\otimes b$, that is
$\nu $, because  $a\otimes b$ is \ {\bf D}-tight. Thus, we get
indeed a function in $DC(\lambda ,\mu ,\nu )$.

Vice versa, suppose we are given a function $f\in DC(\lambda ,\mu
,\nu )$. Then we set  the array $b$ as the array of mixed
derivatives $\partial \partial f$. The array  $a$ is just
diag$(\lambda )$. It is clear that we obtain a standard pair of
type  $(\lambda ,\mu ,\nu )$. It is also clear that the above
constructions are invertible. $\Box$\medskip

We claim that, in the course of this bijection between standard
pairs and discrete concave functions, the associativity bijection
for arrays coincides with the functional associativity bijections
\cite{oct}. For it suffices to clarify the construction of the
second couple of standard pairs in terms of discrete concave
functions (for the first pair it is clear).

For the pair of arrays $(b,c)$, we consider the corresponding
function $\int \!\!\int (b\otimes c)$ on the rectangle  $2n\times
n$. We locate this function on the slope rectangular face  $OABC$
(see Picture 5) of the prism (with $b=(2n,n,n)$) and apply the
octahedron recurrence (with the propagation vector  $(1,0,1)$ and
zero values at the faces $OEA$ and $OEDC$, as we did in Section
2). Due to the functional form of RSK (Theorem 2), we obtain the
function $\int \!\!\int {\bf D}(b\otimes c)$ on top face $EABD$
and the function $\int \!\!\int {\bf L}(b\otimes c)$ on the right
triangle $CDB$.

\unitlength=.80mm \special{em:linewidth 0.4pt} \linethickness{0.4pt}
\begin{picture}(120.00,63.00)
\put(40.00,15.00){\vector(4,-1){54.00}}
\put(40.00,15.00){\vector(0,1){35.00}}
\emline{40.00}{40.00}{1}{60.00}{50.00}{2}
\emline{80.00}{30.00}{3}{40.00}{40.00}{4}
\emline{100.00}{40.00}{5}{60.00}{50.00}{6}
\emline{40.00}{15.00}{7}{60.00}{50.00}{8}
\put(37.00,53.00){\makebox(0,0)[cc]{$z$}}
\put(37.00,13.00){\makebox(0,0)[cc]{$O$}}
\put(37.00,40.00){\makebox(0,0)[cc]{$E$}}
\put(60.00,53.00){\makebox(0,0)[cc]{$A$}}
\put(65.00,9.00){\vector(0,1){24.00}}
\put(65.00,34.00){\vector(2,1){19.00}}
\put(85.00,44.00){\vector(4,-1){24.00}}
\put(65.00,9.00){\vector(3,2){43.00}}
\emline{65.00}{9.00}{9}{85.00}{44.00}{10}
\emline{66.00}{33.00}{11}{110.00}{38.00}{12}
\emline{80.00}{30.00}{13}{90.00}{28.00}{14}
\emline{90.00}{28.00}{15}{90.00}{3.00}{16}
\emline{90.00}{3.00}{17}{110.00}{38.00}{18}
\emline{110.00}{38.00}{19}{90.00}{28.00}{20}
\put(97.00,2.00){\makebox(0,0)[cc]{$x$}}
\put(63.00,6.00){\makebox(0,0)[cc]{$C'$}}
\put(62.00,38.00){\makebox(0,0)[cc]{$D'$}}
\put(87.00,47.00){\makebox(0,0)[cc]{$B'$}}
\put(112.00,40.00){\makebox(0,0)[cc]{$B$}}
\put(90.00,31.00){\makebox(0,0)[cc]{$D$}}
\put(87.00,0.00){\makebox(0,0)[cc]{$C$}}
\put(63.00,22.00){\makebox(0,0)[cc]{$\lambda$}}
\put(73.00,40.00){\makebox(0,0)[cc]{$\mu$}}
\put(99.00,43.00){\makebox(0,0)[cc]{$\nu$}}
\put(77.00,25.00){\makebox(0,0)[cc]{$\sigma$}}
\put(87.00,37.00){\makebox(0,0)[cc]{$\tau$}}
\put(85.00,19.00){\makebox(0,0)[cc]{$\pi$}}
\end{picture}\bigskip

\hfill Picture 5. \hfill\medskip

Thus, we have to recognize the functions which we get on the
triangles (grids) $C'D'B'$, $D'B'B$, $B'BC'$ and $BC'D'$ (these
triangles are the faces of the tetrahedron $C'D'B'B$). Recall that
the functional associativity bijection is obtained via the
octahedron recurrence for this tetrahedron \cite{oct} (with the
same propagation vector and modular flats).\medskip

     1. By Theorem 2, on the {\em first face}  $C'D'B'$ we obtain  $\int \!\!\int $
of the array ${\bf L}b=b$.\medskip

     2. On the {\em second face} $D'B'B$ we get the restriction of the function
$\int \!\!\int (b'\otimes c')$ to this triangle (or the square
$D'B'BD$). Recall that we denoted by $b'\otimes c'$ the standard
pair  ${\bf D}(b\otimes c)$. Thus, by modulo adding the function,
 $\int \!\!\int c'$ is located  on this face.\medskip

     3. On the third (slope) face  $B'BC'$ of the tetrahedron is located
the restriction of function  $\int \!\!\int (b\otimes c)$ to this
triangle. It is easy to see (as in the item 2), that this
restriction is equal to  $\int \!\!\int c$ plus the function of
one variable  $\int $(row sums of $b$). For what follows it is
worth to note that the row sums of  $b$ equals the row sums of
$a\otimes b$ minus the row sums of  $a$, or, equivalently, the
shape of  $\mathbf L(a,b)$ minus the shape of  $a$.\medskip

     4. The only non-trivial function is located on the {\em forth} face
$BC'D'$ of the tetrahedron. Namely, we claim that  {\em the
function on the face  $BC'D'$ is ``almost'' coincided with the
function on the face $BCD$, that is the value of the function at
the $(n+i,i,k)$ is equal to the value of the function $\int
\!\!\int {\bf L}(b\otimes c)$ in the point $(i,k)$}.\medskip

Let us postpone proving this claim, and conclude that we get
finally. On the faces  1 and  3 we have the initial data, the
function $\int \!\!\int b$ on the face 1 and the function  $\int
\!\!\int c$ (modulo adding a function of one variable) on the face
3. At the output faces, that is faces 2 and 4, we have  $\int
\!\!\int c'+\int b'$ and $\int \!\!\int {\bf L}(b\otimes c)$. In
order to get the exact case of  \cite{oct}, we have slightly
modify our functions in order to get discrete concave functions on
the faces 1 and 3. For this, we have to add to the polarized
function on the prism  $OEACDB$ the ``one-dimensional'' function
of the $z$-coordinate (vertical axe) equals  $\int a$, or,
equivalently, $\int \lambda $. Thus, we get discrete concave
function  $\int \lambda +\int \!\!\int b$ on the face 1 with
increments  $\lambda ,\mu $ and  $\sigma $, where $\sigma $
denotes the shape of the array $a\otimes b$. On the face 3, we get
the function  $\int \!\!\int c+(\int l-\int a)+\int a=\int
\!\!\int c+\int l$. Since $l={\bf L}(a\otimes b)$, there holds
$\int l=\int \sigma $. Therefore, on the face 3 is located a
discrete concave function with the increments  $\sigma ,\nu $ and
$\pi $, where  $\pi $ denotes the shape of the array $a\otimes b
\otimes  c$. These are the input functions. The output functions
are: on the face 2 is located the function  $\int \!\!\int c'+\int
b'=\int \!\!\int c'+\int \mu $, since the shape of $b'={\bf D}b$
is   $\mu $. That is a discrete concave function with increments
$\mu ,\nu $ ш $\tau $, where $\tau $ denotes the shape of the
array  $b\otimes c$. On the face 4 is located the function  $\int
\!\!\int {\bf L}(b\otimes c)+\int \lambda $, discrete concave
function corresponding to the standard pair $(a,{\bf L}(b\otimes
c))$, with the increments  $\lambda ,\tau $ ш $\pi $.

Finally, we note that adding the function  $\int \lambda $ of the
vertical variable does not affect on the octahedron recurrence,
and therefore, we get the coincidence of the array associativity
bijection and that is in  \cite{oct} (the functional
form).\medskip

To prove the claim, we observe that after adding the function of
the vertical variable  $\int \lambda $, the polarized function on
the tetrahedron $C'D'DCB$ becomes discrete concave (Corollary  5).
Moreover, it is a constant on each segment in the ground of the
octahedron, which is parallel to the axe  $x$. Therefore, the
function is constant on each segment parallel to the axe  $x$.

Thus, we have prove the following proposition.\medskip

{\bf Theorem 3}. {\em  Under the above isomorphism between  $SP$
and  $DC$, the associativity bijection for array (Section 7)
coincides with the functional associativity bijection constructed
in \cite{oct} using the OR.}

\section{ Commutativity bijection }

Now we turn to the commutativity bijection. Namely, we claim
existence of a natural commutativity bijection between the sets
$SP(\lambda ,\mu ,\nu )$ and $SP(\mu ,\lambda ,\nu )$. Here, as
usual, $n$-tuples  $\lambda ,\mu ,\nu $ denote partitions (not
necessary integer-valued). Thus, we have to associate to a given
standard pair $(a,b)$, of type $(\lambda ,\mu ,\nu )$, a standard
pair $(b',a')$ of type $(\mu ,\lambda ,\nu )$. In \cite{umn}, we
proposed such a bijection, and called it the {\em commuter}.

Here is convenient to use  {\em anti-standard pairs} $(a,b)$, that
is an {\bf R}-tight array $a$ and  an {\bf L}-tight array $b$ such
that  $a\otimes b$ is a {\bf D}-tight array. The type of such an
array is defined similarly as for standard pair, the starting
shape equals the shape of  $a$, the intermediate shape equals the
shape of $b$, and the final shape is the shape of $a\otimes
b$.\medskip

For example, the following concatenated array  $(a,b)$
$$
\begin{array}{|c|c|c||c|c|c|}\hline
 0  & 0 & 0 & 1 & 2 & 1 \\ \hline
  0 & 0 & 2 & 1 & 2 & 0 \\ \hline
  0 & 2 & 1 & 3 & 0 & 0 \\ \hline
\end{array}
$$
is an anti-standard array of the type  $\lambda =(3,2,0)$, $\mu
=(5,4,1)$ ш $\nu =(6,5,4)$.\medskip

The set of anti-standard of type $(\lambda ,\mu ,\nu )$ we let to
denote by $ASP(\lambda ,\mu ,\nu )$. There is a canonical
bijection
$$
SP(\lambda ,\mu ,\nu ) \ \tilde{\to} \ ASP(\lambda ,\mu ,\nu ), \ \
\ (a,b) \mapsto ({\bf R}a,b);
$$
(and the reverse mapping is  $(a,b) \mapsto ({\bf L}a,b)$).

Now, we define the commutativity bijection, commuter,
$$
     Com:ASP(\lambda ,\mu ,\nu ) \to ASP(\mu ,\lambda ,\nu ).
$$
Let $(a,b)$ be an anti-standard pair. We define the commuter from
the rule
$$Com(a,b)={\bf D}(*(a,b))={\bf D}(*b,*a).$$ (Here and in
what follows,  $*$ denotes the central symmetry  of an array,
$*a(i,j)=a(n-i+1,m-j+1)$.) It is easy to check that the pair
$(b',a')={\bf D}(*b,*a)$ is an anti-standard of the required type.
for example, consider $b'={\bf D}(*b)$. Since  $b$ is a {\bf
L}-tight, that is $b={\bf L}b$, we get that  $*b={\bf R}(*b)$ is
{\bf R}-tight. Due to commuting  {\bf D} and  {\bf R}, the array
$b'$ is {\bf R}-tight. Its shape is equal to the shape of  $*b$
and is equal to the shape of  $b$, that is nothing but $\mu $.
Similarly, one can check that we get the correct type.\medskip

{\bf Example}. Let $(a,b)$ be a standard pair from the previous
example. The inverse assay (with respect to  $*$)  $*(a\otimes
b)=*b\otimes *a$ takes the form
$$
\begin{array}{|c|c|c||c|c|c|}\hline
  0 & 0 & 3 & 1 & 2 & 0 \\ \hline
  0 & 2 & 1 & 2 & 0 & 0 \\ \hline
  1 & 2 & 1 & 0 & 0 & 0 \\ \hline
\end{array} \ .
$$
It is easy to check that {\bf D}-condensation of this array equals
$$
\begin{array}{|c|c|c||c|c|c|}\hline
  0 & 0 & 1 & 1 & 2 & 0 \\ \hline
  0 & 1 & 3 & 1 & 0 & 0 \\ \hline
  1 & 3 & 1 & 1 & 0 & 0 \\ \hline
\end{array} \ .
$$
It is easy to check that if we centrally symmetrically reverse the
array and then condensate it, we get the initial array
$$
\begin{array}{|c|c|c||c|c|c|}\hline
 0  & 0 & 0 & 1 & 2 & 1 \\ \hline
  0 & 0 & 2 & 1 & 2 & 0 \\ \hline
  0 & 2 & 1 & 3 & 0 & 0 \\ \hline
\end{array} \ .
$$
Of course, this coincidence has a reason and there holds
\medskip

{\bf Lemma}. {\em  The commuter $Com$ is an involution.}\medskip

That is the composition
$$
             ASP(\lambda ,\mu ,\nu ) \stackrel{Com}{\longrightarrow} ASP(\mu
,\lambda ,\nu ) \stackrel{Com}{\longrightarrow}  ASP(\lambda ,\mu
,\nu )
$$
is the identical mapping.

{\em Proof}. In fact, let $(b',a')=Com(a,b)={\bf D}(*(a\otimes
b))$. Then $Com(b',a')={\bf D}(*(b'\otimes a'))={\bf D}(*{\bf D}
(*(a\otimes b))={\bf D}({\bf U}(a\otimes b))$. According to
\cite{umn}, Section 5.10, the latter equals  ${\bf D}(a\otimes b)$
and this equals to  $a\otimes b$, since $a\otimes b$ is {\bf
D}-tight. \hfill$\Box$\medskip

The corresponding mapping $SP(\lambda,\mu,\nu) \to
SP(\mu,\lambda,\nu)$ we also denote by  $Com$. We claim that, in
integer-valued case, this mapping  $Com$ coincides with the
mapping, which is called the {\em first fundamental symmetry} in
\cite{PV} is denoted by  $\rho_1$. To clarify this claim, we have
to translate the stuff from the array language in the language of
Young tableaux, because the mapping  $\rho_1$ is defined in
\cite{PV} in this language. (Note that this language forces to
restrict ourselves to integer arrays and partitions.) For this
translation we use the following natural bijection
$$
SP_{\mathbb{Z}}(\lambda ,\mu ,\nu ) \to LR(\nu \setminus \lambda
,\mu ),
$$
where  $LR(\nu \setminus \lambda ,\mu )$ denotes the set of
Littlewood-Richardson tableaux of the skew shape  $\nu \setminus
\lambda$ and weight  $\mu$ (\cite{umn,M,PV}).

Let  $(a,b)$ be a standard pair of integer-valued arrays from $
SP_{\mathbb{Z}}(\lambda ,\mu ,\nu )$. Since the array  $a\otimes
b$ is {\bf D}-tight, there is a corresponding semistandard Young
tableau of shape  $\nu $ (see \cite{umn} and end of the item 4 in
Section 5). Moreover,  the {\bf D}-tight array  $a$ defines a
sub-tableau of shape  $\lambda $. The complement to the latter
tableau gives a semistandard skew Young tableau of shape  $\nu
\setminus \lambda $ and weight $\mu $ ($=$shape of $b$), and,
finally, due to  {\bf L}-tightness of the array $b$, this filling
gives a reverse lattice (or dominated or Yamanuchi) word  (see
\cite{umn}, (9.6)). Let us note, that the skew tableau is filled
from the alphabet $n+1,\ldots, 2n$, and, in order to get a
tableau, filled from $1,\ldots, n$, we have to subtract $n$ from
each letter of the tableau.

The reverse mapping is as follows. Let we have a tableau from
$LR(\nu \setminus \lambda ,\mu )$ filled from the alphabet
$1,\ldots, n$. Then we add $n$ to each letter of the tableau and
fill the ``empty part'' of shape $\lambda $ as the Yamanuchi
tableau (i.e. the $j$-th row is constituted only of the letters
$j$, $j=1,\ldots, m$). Thus, we get a semistandard tableau of
shape $\nu$, and the corresponding  {\bf D}-tight array $a\otimes
b$. The array $a$ correspond to the Yamanuchi tableau and,
therefore, is  {\bf L}-tight. Since the word  $w(b)$ is a reverse
lattice word, the array $b$ is {\bf L}-tight.\medskip

For example, for the pair $(a,b)$ from the above example
(specifically, for the standard pair $({\bf L}(a),b)$) we obtain
the following tableau

\unitlength=.5mm \special{em:linewidth 0.4pt} \linethickness{0.4pt}
\begin{picture}(90.00,45.00)(-55,0)
\emline{30.00}{5.00}{1}{30.00}{35.00}{2}
\emline{30.00}{35.00}{3}{70.00}{35.00}{4}
\emline{70.00}{35.00}{5}{70.00}{25.00}{6}
\emline{70.00}{25.00}{7}{80.00}{25.00}{8}
\emline{80.00}{25.00}{9}{80.00}{15.00}{10}
\emline{80.00}{15.00}{11}{90.00}{15.00}{12}
\emline{90.00}{15.00}{13}{90.00}{5.00}{14}
\emline{90.00}{5.00}{15}{30.00}{5.00}{16}
\put(35.00,30.00){\makebox(0,0)[cc]{4}}
\put(45.00,30.00){\makebox(0,0)[cc]{5}}
\put(55.00,30.00){\makebox(0,0)[cc]{5}}
\put(65.00,30.00){\makebox(0,0)[cc]{6}}
\put(35.00,20.00){\makebox(0,0)[cc]{2}}
\put(45.00,20.00){\makebox(0,0)[cc]{2}}
\put(55.00,20.00){\makebox(0,0)[cc]{4}}
\put(65.00,20.00){\makebox(0,0)[cc]{5}}
\put(75.00,20.00){\makebox(0,0)[cc]{5}}
\put(35.00,10.00){\makebox(0,0)[cc]{1}}
\put(45.00,10.00){\makebox(0,0)[cc]{1}}
\put(55.00,10.00){\makebox(0,0)[cc]{1}}
\put(65.00,10.00){\makebox(0,0)[cc]{4}}
\put(75.00,10.00){\makebox(0,0)[cc]{4}}
\put(85.00,10.00){\makebox(0,0)[cc]{4}}
\emline{60.00}{5.00}{17}{60.00}{15.00}{18}
\emline{60.00}{15.00}{19}{50.00}{15.00}{20}
\emline{50.00}{15.00}{21}{50.00}{25.00}{22}
\emline{50.00}{25.00}{23}{30.00}{25.00}{24}
\end{picture}

\noindent One can check that the word $w(b)=4556455444$ is a
Yamanuchi word indeed (in the alphabet $\{4,5,6\}$).\medskip

Recall that in \cite{PV} the first fundamental  symmetry is
defined using the tableaux-switching algorithm. In this example,
we have to transport the letters  1, 2 and  3 through the letters
4, 5, 6. As a result, we get a tableau  with the letters order:
$4<5<6<1<2<3$).\medskip

\unitlength=.500mm \special{em:linewidth 0.4pt}
\linethickness{0.4pt}
\begin{picture}(90.00,35.00)(-40,0)
\emline{30.00}{5.00}{1}{30.00}{35.00}{2}
\emline{30.00}{35.00}{3}{70.00}{35.00}{4}
\emline{70.00}{35.00}{5}{70.00}{25.00}{6}
\emline{70.00}{25.00}{7}{80.00}{25.00}{8}
\emline{80.00}{25.00}{9}{80.00}{15.00}{10}
\emline{80.00}{15.00}{11}{90.00}{15.00}{12}
\emline{90.00}{15.00}{13}{90.00}{5.00}{14}
\emline{90.00}{5.00}{15}{30.00}{5.00}{16}
\emline{40.00}{35.00}{17}{40.00}{25.00}{18}
\emline{40.00}{25.00}{19}{70.00}{25.00}{20}
\emline{70.00}{25.00}{21}{70.00}{15.00}{22}
\emline{70.00}{15.00}{23}{80.00}{15.00}{24}
\emline{80.00}{15.00}{25}{80.00}{5.00}{26}
\put(35.00,30.00){\makebox(0,0)[cc]{6}}
\put(45.00,30.00){\makebox(0,0)[cc]{1}}
\put(55.00,30.00){\makebox(0,0)[cc]{2}}
\put(65.00,30.00){\makebox(0,0)[cc]{2}}
\put(35.00,20.00){\makebox(0,0)[cc]{5}}
\put(45.00,20.00){\makebox(0,0)[cc]{5}}
\put(55.00,20.00){\makebox(0,0)[cc]{5}}
\put(65.00,20.00){\makebox(0,0)[cc]{5}}
\put(75.00,20.00){\makebox(0,0)[cc]{1}}
\put(35.00,10.00){\makebox(0,0)[cc]{4}}
\put(45.00,10.00){\makebox(0,0)[cc]{4}}
\put(55.00,10.00){\makebox(0,0)[cc]{4}}
\put(65.00,10.00){\makebox(0,0)[cc]{4}}
\put(75.00,10.00){\makebox(0,0)[cc]{4}}
\put(85.00,10.00){\makebox(0,0)[cc]{1}}
\end{picture}

This tableau corresponds to the pair  $(b',a')$ of the
above-mentioned example, and this is not a curious. Specifically,
we claim that there holds the following
\medskip

{\bf Proposition 3}.  {\em The above defined bijection between
$SP_{\mathbb{Z}}(\lambda ,\mu ,\nu )$ and  $LR(\nu \setminus
\lambda ,\mu )$ makes the following diagram commutative (that is
coincidence  $Com$ and  $\rho_1$) }
$$
\begin{array}{ccc}
{SP_{\mathbb{Z}}(\lambda ,\mu ,\nu )}&{\rightarrow} & {LR(\nu
\setminus \lambda ,\mu )}\cr {Com\downarrow} & &
{{\rho_1}\downarrow}\cr {SP_{\mathbb{Z}}(\mu,\lambda ,\nu
)}&{\rightarrow} & {LR(\nu \setminus \mu, \lambda )}\cr
\end{array}
$$\medskip

In fact, the tableaux-switching might be obtained (\cite{PV}, 3.1)
of the form of the composition of three Schutzenberger involutions
${\bf S}_1 {\bf S}_{12} {\bf S}_1$. That is, for a standard pair
$(a,b)$, we, first, apply  ${\bf S}$  to the array $a$, then we
apply  ${\bf S}$ to the array  ${\bf S}a\otimes b$, and finally,
we apply  ${\bf S}$ to the first $n\times n$ part of the array
${\mathbf S}({\bf S}a\otimes b)$.

According to \cite{umn}, the array  ${\bf S}a$ is equal to  ${\bf
D}(*a)$. Since  $a$ is bi-tight array ({\bf D}- and {\bf
L}-tight), the array  $*a$ is {\bf R}-tight. Therefore ${\bf
D}(*a)$ is a (unique) {\bf D}- and {\bf R}-tight array of shape
$\lambda $. But the array ${\bf R}a$ is also {\bf D}- and {\bf
R}-tight of shape $\lambda $. Thus, we get ${\bf S}a={\bf
D}(*a)={\bf R}a$, and the application of ${\mathbf S}_1$ is
nothing but sending a standard pair to the antistandard one.

Now,  ${\bf S}={\bf S}_{12}({\bf R}a\otimes b)$ is equal to  ${\bf
D}(*({\bf R}a\otimes b))$, and finally, applying ${\bf S}_1$ sends
the antistandard pair  ${\bf D}(*({\bf R}a\otimes b))$ to a
standard one. This composition does exactly that $Com$ does.
$\Box$

\section{Functional form of the commutativity bijection}

Now we translate the commutativity bijection into the language of
functions. Recall (see Section 8) that the set $SP(\lambda ,\mu
,\nu )$ (and also $ASP(\lambda ,\mu ,\nu )$) is bijective to the
set $DC(\lambda ,\mu ,\nu )$ of discrete concave functions on the
triangle grid with the increments $\lambda ,\mu ,\nu $ (increments
along the left-hand side constitute an $n$-tuple $\lambda $, along
the top of the triangle constitute $\mu $, and $\nu$ along the
hypotenuse).

Using these bijections and the bijection
$$
 Com:ASP(\lambda ,\mu ,\nu ) \to ASP(\mu ,\lambda ,\nu ),
$$
we obtain (from the commutative diagram) a bijection $Com'$
between the set  $DC(\lambda ,\mu ,\nu )$ and $DC(\mu , \lambda
,\nu )$. Here we transform this definition of $Com'$ in more
direct and transparent form.

Let $f\in DC(\lambda,\mu,\nu)$. Consider the corresponding  {\bf
L}-tight array $b=\partial \partial f$; the vector of its column
sums ($I$-weight) is equal to  $\mu $. The vector of row sums
($J$-weight) is equal to $\nu -\lambda $. Let us pick the  {\bf
D}{\bf R}-tight array  $a={\bf R}(\text{diag}(\lambda))$ of shape
$\lambda $. The pair  $(a,b)$ is anti-standard and corresponds to
$f$.\medskip

The commuter  $Com$ sends the anti-standard pair $(a,b)$ to the
anti-standard pair  $(b',a')={\bf D}(*b,*a)$. (Of course, we know
that  $b'$ is the {\bf D}{\bf R}-tight array of shape  $\mu $, and
we know that the vector of column sums of $a'$ is equal to that of
$*a$, and that is indeed $\lambda$. In fact, the vector of column
sums of $a$ is equal to $\lambda^{op}=(\lambda_n,...,\lambda_1)$,
therefore, that of  $*a$ is equal to  $\lambda $.) Now, from this
pair we have to return to a discrete concave function, and, by the
definition,  we get
$$
Com'(f)=\int \! \! \!\!\int a'+\int \mu \ .
$$

To understand better the pair $(b',a')$ we exploit Theorem 2. Pick
the supermodular function  $g=\int\!\!\int *b\otimes*a$, and
locate it on the slope face $OABC$ (see Picture  6), and then
apply the octahedron recurrence in the prism as in Section 6. Then
on the top face $EABD$ we get the function which corresponds to
$(b',a')$. Specifically, the function  $Com'(f)$ is located on the
triangle $D'B'B$ (the increments along the left side $D'B'$
constitute  $\mu$, $\lambda$ along  the top side $B'B$, and $\nu$
along the hypotenuse $D'B$).

\unitlength=.600mm \special{em:linewidth 0.4pt}
\linethickness{0.4pt}
\begin{picture}(120.00,90.00)(-50,0)
\put(20.00,30.00){\vector(3,-1){75.00}}
\put(20.00,30.00){\vector(0,1){45.00}}
\put(24.00,76.00){\makebox(0,0)[cc]{$z$}}
\put(98.00,8.00){\makebox(0,0)[cc]{$x$}}
\emline{20.00}{70.00}{1}{50.00}{80.00}{2}
\emline{50.00}{80.00}{3}{110.00}{60.00}{4}
\emline{110.00}{60.00}{5}{80.00}{10.00}{6}
\emline{80.00}{10.00}{7}{80.00}{50.00}{8}
\emline{80.00}{50.00}{9}{20.00}{70.00}{10}
\emline{50.00}{20.00}{11}{50.00}{60.00}{12}
\emline{50.00}{60.00}{13}{80.00}{70.00}{14}
\emline{80.00}{70.00}{15}{50.00}{20.00}{16}
\emline{50.00}{20.00}{17}{110.00}{60.00}{18}
\emline{110.00}{60.00}{19}{80.00}{50.00}{20}
\put(50.00,22.00){\vector(0,1){36.00}}
\put(53.00,61.00){\vector(3,1){24.00}}
\put(53.00,22.00){\vector(3,2){52.00}}
\put(83.00,69.00){\vector(3,-1){24.00}}
\put(17.00,27.00){\makebox(0,0)[cc]{$O$}}
\put(48.00,17.00){\makebox(0,0)[cc]{$C'$}}
\put(77.00,7.00){\makebox(0,0)[cc]{$C$}}
\put(17.00,71.00){\makebox(0,0)[cc]{$E$}}
\put(50.00,64.00){\makebox(0,0)[cc]{$D'$}}
\put(80.00,53.00){\makebox(0,0)[cc]{$D$}}
\put(52.00,82.00){\makebox(0,0)[cc]{$A$}}
\put(82.00,72.00){\makebox(0,0)[cc]{$B'$}}
\put(112.00,63.00){\makebox(0,0)[cc]{$B$}}
\emline{20.00}{30.00}{21}{25.00}{38.00}{22}
\emline{30.00}{46.00}{23}{35.00}{54.00}{24}
\emline{40.00}{62.00}{25}{45.00}{70.00}{26}
\emline{50.00}{80.00}{27}{47.00}{74.00}{28}
\emline{50.00}{60.00}{29}{110.00}{60.00}{30}
\put(97.00,73.00){\vector(-3,-2){15.00}}
\put(102.00,76.00){\makebox(0,0)[cc]{$Com'(f)$}}
\end{picture}

\hfill ашё. 6 \hfill\medskip

Now we show that the function $g=\int\!\!\int (*b\otimes *a)$ is
very simple related to the function  $f$. Specifically, we explain
a relation with the function $\int\!\!\int (a\otimes b)$ on the
rectangle $2n\times n$. Since the restriction of this function to
the triangle $C'B'B$ is exactly  $f$, we also denote the function
$\int\!\!\int (a\otimes b)$ by $f$. Its boundary increments we
depicted on Picture 7.\medskip

\unitlength=1mm \special{em:linewidth 0.4pt} \linethickness{0.4pt}
\begin{picture}(83.00,30.00)
\emline{40.00}{5.00}{1}{80.00}{5.00}{2}
\emline{80.00}{5.00}{3}{80.00}{4.00}{4}
\put(80.00,5.00){\vector(0,1){20.00}}
\put(40.00,5.00){\vector(0,1){20.00}}
\put(40.00,25.00){\vector(0,0){0.00}}
\put(40.00,25.00){\vector(1,0){20.00}}
\put(60.00,25.00){\vector(0,0){0.00}}
\put(60.00,25.00){\vector(1,0){20.00}}
\emline{60.00}{5.00}{5}{60.00}{25.00}{6}
\emline{60.00}{5.00}{7}{80.00}{25.00}{8}
\put(37.00,15.00){\makebox(0,0)[cc]{$0$}}
\put(50.00,2.00){\makebox(0,0)[cc]{$0$}}
\put(70.00,2.00){\makebox(0,0)[cc]{$0$}}
\put(50.00,28.00){\makebox(0,0)[cc]{$\lambda^{op}$}}
\put(70.00,28.00){\makebox(0,0)[cc]{$\mu$}}
\put(83.00,15.00){\makebox(0,0)[cc]{$\nu$}}
\put(67.00,17.00){\makebox(0,0)[cc]{$f$}}
\end{picture}

\hfill Picture 7. Function  $f$. \hfill \medskip

If we ``turn over'' the function $f$, that is if we consider the
function  $*f$, given by the rule $(*f)(i,j)=f(2n-i,n-j)$, then we
will almost get the function  $g$. More precisely, $*f$ and $g$
have the same mixed derivatives $\partial \partial $, but are
differ in the boundary increments and have different values at
$(0,0)$ (see the next picture).
\medskip

\unitlength=1mm \special{em:linewidth 0.4pt} \linethickness{0.4pt}
\begin{picture}(83.00,30.00)
\put(40.00,5.00){\vector(1,0){20.00}}
\put(60.00,5.00){\vector(1,0){20.00}}
\put(40.00,5.00){\vector(0,1){20.00}}
\emline{40.00}{25.00}{1}{80.00}{25.00}{2}
\emline{80.00}{25.00}{3}{80.00}{5.00}{4}
\emline{60.00}{5.00}{5}{60.00}{25.00}{6}
\emline{40.00}{5.00}{7}{60.00}{25.00}{8}
\put(54.00,12.00){\makebox(0,0)[cc]{$*f$}}
\put(35.00,15.00){\makebox(0,0)[cc]{$-\nu^{op}$}}
\put(50.00,28.00){\makebox(0,0)[cc]{$0$}}
\put(70.00,28.00){\makebox(0,0)[cc]{$0$}}
\put(83.00,15.00){\makebox(0,0)[cc]{$0$}}
\put(50.00,2.00){\makebox(0,0)[cc]{$-\mu^{op}$}}
\put(70.00,2.00){\makebox(0,0)[cc]{$-\lambda$}}
\end{picture}

\hfill Picture 8. Function $*f$. \hfill \medskip

Changing of boundary values (with preserving mixed derivatives)
can be make by adding an appropriate separable function (of
variables  $x$ and $z=y$). Doing this, we obtain
$$
 g=*f+\int_x(\mu^{op},\lambda )+\int_z \nu^{op}-f(0,0)=
 *f+\int_x(\mu^{op},\lambda )+\int_z \nu^{op}-|\nu|.
$$
Here  $(\mu^{op},\lambda)$ denotes the tuple
$(\mu_n,...,\mu_1,\lambda_1,...,\lambda_n)$.

\unitlength=1mm \special{em:linewidth 0.4pt} \linethickness{0.4pt}
\begin{picture}(83.00,33.00)
\emline{40.00}{25.00}{1}{40.00}{5.00}{2}
\emline{40.00}{5.00}{3}{80.00}{5.00}{4}
\put(80.00,5.00){\vector(0,1){20.00}}
\put(40.00,25.00){\vector(1,0){20.00}}
\put(60.00,25.00){\vector(1,0){20.00}}
\emline{60.00}{25.00}{5}{60.00}{5.00}{6}
\emline{60.00}{5.00}{7}{80.00}{25.00}{8}
\put(37.00,15.00){\makebox(0,0)[cc]{$0$}}

\put(55.00,15.00){\makebox(0,0)[cc]{$g$}}

\put(50.00,2.00){\makebox(0,0)[cc]{$0$}}
\put(70.00,2.00){\makebox(0,0)[cc]{$0$}}
\put(50.00,28.00){\makebox(0,0)[cc]{$\mu^{op}$}}
\put(70.00,28.00){\makebox(0,0)[cc]{$\lambda$}}
\put(84.00,15.00){\makebox(0,0)[cc]{$\nu^{op}$}}
\emline{80.00}{10.00}{9}{65.00}{10.00}{10}
\put(74.00,12.00){\makebox(0,0)[cc]{$const$}}
\end{picture}

\hfill Picture 9. Function $g$. \hfill \medskip

For what follows it is worth to note that the function $g$ is
constant on the segments $[(n+j,j),(2n,j)]$, because the array
$*a$ takes zero values at these segments.

Thus, let us conclude this part with recalling what we have done.
We locate the function $g$ on the face  $OABC$ of the prism, using
the OR, we obtain a polarized function $G$ on the prism. The
restriction of $G$ to the triangle $D'B'B$ is the function of our
interest $Com'(f)$. Now, we are interested in the restriction of
the function $G$ to the vertical cross connected wall $C'D'B'$. By
Theorem 2, we know that this function is $\int\!\!\int {\bf
L}(*b)$ with boundary increments  $0$, $\mu $ and  $\nu
^{op}-\lambda ^{op}=(\nu -\lambda ) ^{op}$. Let us add to this
function the ``one-dimensional'' function  $\int _z (-\nu^{op})$.
Thus, the function
$$
h=\int\!\!\!\!\int {\bf L}(*b)-\int _z \nu^{op}+|\nu|.
$$
is discrete concave and has increments  $-\nu^{op} , \mu $ and
$-\lambda ^{op}$.\medskip

\unitlength=1mm \special{em:linewidth 0.4pt} \linethickness{0.4pt}
\begin{picture}(75.00,30.00)(30,0)
\put(50.00,5.00){\vector(0,1){20.00}}
\put(50.00,25.00){\vector(1,0){20.00}}
\put(50.00,5.00){\vector(1,1){20.00}}
\put(56.00,17.00){\makebox(0,0)[cc]{$h$}}
\put(48.00,3.00){\makebox(0,0)[cc]{$C'$}}
\put(48.00,27.00){\makebox(0,0)[cc]{$D'$}}
\put(72.00,27.00){\makebox(0,0)[cc]{$B'$}}
\put(45.00,15.00){\makebox(0,0)[cc]{$-\nu^{op}$}}
\put(60.00,27.00){\makebox(0,0)[cc]{$\mu$}}
\put(66.00,15.00){\makebox(0,0)[cc]{$-\lambda^{op}$}}
\put(85.00,15.00){\makebox(0,0)[cc]{\text{or}}}
\end{picture}
\begin{picture}(75.00,30.00)(50,0)
\put(50.00,25.00){\vector(0,-1){20.00}}
\put(50.00,25.00){\vector(1,0){20.00}}
\put(70.00,25.00){\vector(-1,-1){20.00}}
\put(56.00,17.00){\makebox(0,0)[cc]{$h$}}
\put(48.00,3.00){\makebox(0,0)[cc]{$C'$}}
\put(48.00,27.00){\makebox(0,0)[cc]{$D'$}}
\put(72.00,27.00){\makebox(0,0)[cc]{$B'$}}
\put(47.00,15.00){\makebox(0,0)[cc]{$\nu$}}
\put(60.00,27.00){\makebox(0,0)[cc]{$\mu$}}
\put(66.00,15.00){\makebox(0,0)[cc]{$\lambda$}}
\end{picture}

\hfill  Picture 10. Function  $ h$. \hfill \medskip

The inter-relation between the function $h$ and the function
$Com'(f)$ is established in the following proposition.\medskip

{\bf Proposition 4}.  {\em  1) The octahedron recurrence on the
quarter of the octahedron $OD'B'C'$ (the Henriques-Kamnitzer
construction) sends the function  $*f$ (on the face $OB'C'$) to
the function $h$, that is  $h=HK(*f)$ up to adding a constant.

     2) $h$ almost coincides with  $Com'(f)$. Specifically, if we consider
$Com'(f)$ as a function on the triangle  $D'B'B$, then there holds
$h(i,j)=Com'(f)(n-j,n-j+i,n)$.}\medskip

In other words, in order to get $Com'(f)$, we have to rotate the
function $h$ ``counter-clockwise on $120^{\circ}$ '', see Picture
11.

\unitlength=1.00mm \special{em:linewidth 0.4pt}
\linethickness{0.4pt}
\begin{picture}(100.00,32.00)
\put(30.00,25.00){\vector(0,-1){20.00}}
\put(30.00,25.00){\vector(1,0){20.00}}
\put(50.00,25.00){\vector(-1,-1){20.00}}
\put(36.00,17.00){\makebox(0,0)[cc]{$h$}}
\put(40.00,28.00){\makebox(0,0)[cc]{$\mu$}}
\put(26.00,15.00){\makebox(0,0)[cc]{$\nu$}}
\put(45.00,15.00){\makebox(0,0)[cc]{$\lambda$}}
\put(80.00,25.00){\vector(1,0){20.00}}
\put(80.00,5.00){\vector(0,1){20.00}}
\put(80.00,5.00){\vector(1,1){20.00}}
\put(88.00,21.00){\makebox(0,0)[cc]{$Com'(f)$}}
\put(90.00,28.00){\makebox(0,0)[cc]{$\lambda$}}
\put(76.00,15.00){\makebox(0,0)[cc]{$\mu$}}
\put(95.00,15.00){\makebox(0,0)[cc]{$\nu$}}
\put(60.00,12.50){\oval(10.00,7.00)[b]}
\put(65.00,12.00){\vector(0,1){3.00}}
\end{picture}

\hfill Picture 11.  \hfill \medskip

{\bf  Corollary}. {\em  In the language of discrete concave
functions the commuter  $Com'$ coincides with the OR commuter due
to Henriques and Kamnitzer.}\medskip

In remark 4.5 in \cite{h-k} is noted, that in order to get the
HK-commuter, one has to apply the OR on the quarter of the
octahedron, and then to rotate the picture counter-clockwise on
$120^{\circ}$ as in Picture  11.\medskip

{\small Let us illustrate this one the example from the beginning
of this section. It is easy to see, that the corresponding
function  $g$ (which we locate on the slope face of the prism)
takes the form
$$
\begin{array}{ccccccc}
  0 & 1 & 5 & 10 & 13 & 15 & 15 \\
  0 & 1 & 5 & 7 & 9 & 9 & 9 \\
  0 & 1 & 3 & 4 & 4 & 4 & 4 \\
  0 & 0 & 0 & 0 & 0 & 0 & 0
\end{array}
$$
The restriction of the function $G$ to the tope face is equal to
the following function
$$
\begin{array}{ccccccc}
  0 & 1 & 5 & 10 & 13 & 15 & 15 \\
  0 & 1 & 5 & 9 & 11 & 11 & 11 \\
  0 & 1 & 4 & 5 & 6 & 6 & 6 \\
  0 & 0 & 0 & 0 & 0 & 0 & 0
\end{array} ,
$$
and, in particular, we get  $Com'(f)$
$$
\begin{array}{cccc}
  10 & 13 & 15 & 15 \\
  9 & 11 & 11 &  \\
  5 & 6 &  &  \\
  0 &  &  &
\end{array} .
$$
The following function is equal to the restriction of $G$ to the
triangle $C'D'B'$:
$$
\begin{array}{cccc}
  0 & 5 & 9 & 10 \\
  0 & 5 & 7 &  \\
  0 & 4 &  &  \\
  0 &  &  &
\end{array} .
$$
And after adding to this function $z$-function $-\int \nu +|\nu|$,
we get the function $h$
$$
\begin{array}{cccc}
  0 & 5 & 9 & 10 \\
  6 & 11 & 13 &  \\
  11 & 15 &  &  \\
  15 &  &  &
\end{array} .
$$
From these computations one can see that  $h$ coincides with
$Com'(f)$ after $120^{\circ}$ counter-clockwise rotation.}\medskip

{\em   Proof of Proposition 4}. 1) By Theorem 2, we get the
function  $\int\!\!\int {\bf L}(*b)$ (on $C'D'B'$) from the
restriction of  $g$ to $OAB'C'$ and the OR. Since the OR commutes
with adding (and subtracting) of separable functions of variables
$x$ and $z$), we subtract from  $g$ the separable function $\int_x
(\mu ^{op} ,\lambda )+\int_z \nu ^{op}-|\nu|$. Then on $OAB'C'$ we
get the function $*f$. And we get the function $h$ (up to
constant) on the wall $C'D'B'$, since we have to subtract $\int _z
\nu^{op}$  from $\int\!\!\int {\bf L}(*b)$.

Thus on the quarter of the octahedron, the octahedron recurrence
sends $*f$ into $h$.\medskip

     2) Adding a function of the $z$-variable to the function $G$
does not affect on the top face $EABD$. Therefore, we add the
function  $-\int _z \nu ^{op}+|\nu|$ to the function $G$.

Let us consider the restriction of the function  $G-\int _z \nu
^{op}+|\nu|$ to the tetrahedron $C'D'B'B$. It is clear that this
is a polarized function. Thereupon, its restriction to the face
$C'D'B'$, being the function $\int\!\!\int {\bf L}(*b)-\int _z\nu
^{op}+|\nu|=h$, is a discrete concave function. In fact, by modulo
of adding an affine function, this function corresponds to the
standard pair  $(\mbox{diag}(-\nu^{op}) ,{\bf L}(*b))$, and, thus,
is discrete concave. Finally, the restriction to the face  $C'B'B$
is also discrete concave, and, moreover, it is constantly equal
$|\nu|$ on the edge  $C'B$. Let us verify this claim.

The restriction of  $g$ to the square  $C'B'BC$ is equal to the
function  $\int\!\!\int *a$ plus the restriction of $g$ to the
edge $C'B'$, that is  $\int_z (\nu ^{op}-\lambda ^{op})$. Of our
interest is the function  $g-\int_z \nu^{op}$, that is the
function $\int\!\!\int *a-\int_z\lambda ^{op}$. Now we have
honestly deduce the claim from the $\mathbf {LU}$-tightness of
$*a$ ({\bf DR}-tightness of $a$).

Note that it suffices to check the claim for  $\lambda$ of the
form $(1,...,1,0,...,0)=(1^k,0^{n-k})$. In this case, one can find
values of the function $\int\!\!\int *a-\int_z\lambda ^{op}$ at a
point $(i,j)$, that is  $\min(k+i-j,0)$. From this discrete
concavity is obvious. On the edge $C'B$, we have $i=j$, and our
function is equal to  $\min(k,0)=0$, since  $k\ge 0$.

Thus, the polarized function $G-\int _z \nu^{op}+|\nu|$ has
discrete concave restrictions to the faces $C'D'B'$ and  $C'B'B$
of the tetrahedron $C'D'B'B$. By Theorem 1 $G$ is polarized
discrete concave function. Since $G$ is a constant on the edge
$C'Т$, $G$ is a constant function on any segment, parallel to
$(1,1,1)$, of the tetrahedron. In particular, the values  in the
points $(n,i,j)$ and $(n,i,j)+(n-j)(1,1,1)=(2n-j,n+i-j,n)$
coincide. But we claimed exactly that in the item 2. $\Box$

\section{Commuter and the second fundamental symmetry}

Here we prove that, for the integer-valued set-up, the commuter
$Com$ coincides with the second fundamental symmetry due to Pak
and Vallexo. This bijection (we will consider  $\rho_2'$) is
defined as the following composition (for details see \cite{PV})
$$
      LR(\nu \setminus \lambda ,\mu ) \stackrel{\tau}{\to} T(\mu ,\nu
-\lambda ) \stackrel{\bf S}{\to}T(\mu ,(\nu -\lambda )^{op}  )
\stackrel{\gamma^{-1}}{\to} LR(\nu \setminus \mu ,\lambda ).
$$
Here $T(\mu ,\alpha)$ denotes some subset of Young tableaux of the
shape $\mu $ and the weight $\alpha$.

Let us translate this composition into the language of arrays.
Recall, that a tableau of  $LR(\nu \setminus \lambda ,\mu )$ might
be seen as an  {\bf L}-tight array $b$, such that (diag$(\lambda),
b)$ is a standard pair (that is diag$(\lambda)\otimes b$ is {\bf
D}-tight). The mapping  $\tau$ makes transposition of  $b$, and
gives a  {\bf D}-tight array $b^T$. The Schutzenberger involution
$\bf S$ sends $b^T$ into ${\bf D}(*b^T)$. After repeated
transposition, we obtain {\bf L}-tight array ${\bf D}(*b^T )^T
={\bf L}(*b^{TT})={\bf L}(*b)$. Note, that we have yet met this
array in the previous Section! Recall that the array diag$(-\nu
^{op} )\otimes {\bf L}(*b)$ is $\mathbf D$-tight and $h$ is equal
to $\int \! \!\int {\bf L}(*b)-\int \nu^{op}+|\nu|$.

The commuter $a'=Com(b)$ gives the function $\tilde
h=Com'(f)=\int\! \!\int a'+\int \mu$. By Proposition 4, this
function is easily recalculated from the function  $h$: $\tilde
h(n-j,n-j+i)=h(i,j)$. Thus, it remains to explain a relation
between the arrays $a'=\partial\partial \tilde h$ and $\tilde
b=Com(b)=\partial \partial h$. Specifically, we express $\tilde b$
by means of $a'$. By the definition
$$
\tilde b(i,j)=h(i,j)-h(i-1,j)-h(i,j-1)+h(i-1,j-1).
$$
Recalling the relation between  $h$ and $\tilde h$, we get
$$
\tilde b(i,j)= $$ $$\tilde h(n-j,n-j+i)-\tilde h(n-j,n-j+i+1)-\tilde
h(n-j+1,n-j+i+1)+\tilde h(n-j+1,n-j+i).
$$
The first two summands in the above expression give
$$
\mu_{n-j+i}+a'(1,n-j+i)+a'(2,n-j+i)+...+a'(n-j,n-j_i).
$$
Analogously, one can find the difference of the third and the
forth summands. Thus, we get
$$
\tilde
b(i,j)=[\mu_{n-j+i}+a'(1,n-j+i)+a'(2,n-j+i)+...+a'(n-j,n-j+i)]-
$$
$$
[\mu_{n-j+i+1}+a'(1,n-j+i+1)+a'(2,n-j+i+1)+...+a'(n-j+1,n-j+i+1)].
$$
This is exactly the definition of the mapping $\gamma$ in
\cite{PV}. Thus, we get that, under bijection between
$LR(\nu\setminus \lambda,\mu)$  and $SP_{\mathbb
Z}(\lambda,\mu,\nu)$, the commuter  $Com$  coincides with the
mapping  $\rho_2'$. Due to Lemma the commuter $Com$ is an
involution, that is it coincides with its reversion. In \cite{PV}
is shown that $\rho_2'$ is the reversion to  $\rho_2$. Thus
$\rho_2=\rho_2'$ and we have proved the following theorem.\medskip

{\bf Theorem 4}. {\em  The commuter $Com$ (defined in terms of
arrays) coincides with the Henriques-Kamnitzer commuter (defined
in terms of discrete concave functions), and coincides with the
Pak-Vallexo fundamental symmetries  $\rho _1 $, $\rho _2$ and
$\rho _2 '$ (defined in terms of Young tableaux).}\medskip

V.I.Danilov: CEMI RAS, Nahimovskii pr. 47, 117418 Moscow, Russia,
e-mail: vdanilov43@mail.ru\bigskip

G.A. Koshevoy: CEMI RAS and Poncelet Laboratory (UMI 2615 of CNRS
and Independent University of Moscow),  e-mail:
koshevoy@cemi.rssi.ru

\end{document}